\documentclass{amsart}

\usepackage{pb-diagram}
\usepackage{pstricks}
 \newtheorem{thm}{Theorem}[section]
 
 \newtheorem{lem}[thm]{Lemma}
 \newtheorem{prop}[thm]{Proposition}
 
 \theoremstyle{definition}
 \newtheorem{defn}[thm]{Definition}
 \theoremstyle{definition}
 
 \theoremstyle{remark}
 \newtheorem{rem}[thm]{Remark}
 \numberwithin{equation}{section}
 \newtheorem{expl}[thm]{Example}

\begin{document}

\title[Infinity-Inner-Products on A-Infinity-Algebras]
{Infinity-Inner-Products on A-Infinity-Algebras}
\author{Thomas Tradler}

\begin{abstract}
In this paper the Hochschild-cochain-complex of an
A$_\infty$-algebra $A$ with values in an A$_\infty$-bimodule $M$
over $A$ and maps between them is defined. Then, an
$\infty$-inner-product on $A$ is defined to be an
A$_\infty$-bimodule-map between $A$ and its dual $A^{*}$.
There is a graph-complex associated to A$_\infty$-algebras with
$\infty$-inner-product.
\end{abstract}

\maketitle

\tableofcontents

\section{Introduction}

An A$_\infty$-algebra $(A,D)$ consists of a module $A$ (over a
given ring $R$) together with a coderivation $D$ on the shifted
tensor-coalgebra
$$T(sA):=\bigoplus_{i\geq 0} (sA)^{\otimes i},$$
$$\Delta(sa_{1},...,sa_{n}):=\sum_{i=0}^{n} (sa_{1},...,sa_{i})
    \otimes(sa_{i+1},...,sa_{n}),$$
such that $D^{2}=0$. The Hochschild-cochain-complex of $A$ is
given by the space $CoDer(T(sA),T(sA))$ of coderivations on
$T(sA)$, together with the differential $\delta(f):= [D,f]=D\circ
f\pm f\circ D$. This will be reviewed in section 2.

In section 3, the notion of an A$_\infty$-bimodule $(M,D^{M})$
over the A$_\infty$-algebra $(A,D)$ is defined, by taking a
coderivation $D^{M}$ from the coalgebra $T(sA)$ into the
bi-comodule $T^{sM}(sA)$, which is given by
$$T^{sM}(sA):=R\oplus\bigoplus_{k\geq 0, l\geq 0} (sA)^{\otimes k}
\otimes (sM) \otimes (sA)^{\otimes l},$$
\begin{eqnarray*}
\Delta^{sM}(sa_{1},...,sm,...,sa_{k+l})&:=&
 \sum_{i=0}^{k}(sa_{1},...,sa_{i}) \otimes(sa_{i+1},...,sm,...,sa_{n})+\\
 & &+\sum_{i=k}^{k+l} (sa_{1},...,sm,...,sa_{i})
 \otimes(sa_{i+1},...,sa_{k+l}),
\end{eqnarray*}
such that $(D^{M})^{2}=0$.
Then this definition of $D^{M}$ implies that the differential
$\delta^{M}$ on $CoDer(TsA,T^{sM}sA)$, given by $\delta^{M}(f)
:=D^{M}\circ f\pm f\circ D$,  is well-defined. So, one can call
$CoDer(TA,T^{M}A)$ together with the differential $\delta^{M}$ the
Hochschild-cochain-complex of $A$ with values in the bialgebra
$M$.

In section 4, the concept of an A$_\infty$-bimodule-map between
two A$_\infty$-bimodules $(M,D^{M})$ and $(N,D^{N})$ over $(A,D)$
is defined. This is given by a map $F: T^{sM}sA\longrightarrow
T^{sN}sA$ that respects the differentials: $D^{N}\circ F=F\circ
D^{M}$. Again $F$ has to have good properties with respect to the
coderivations $\Delta^{M}$ and $\Delta^{N}$. These properties will
make the induced map $F^{\sharp}:CoDer(TsA,T^{sM}sA)
\longrightarrow CoDer(TsA, T^{sN}sA)$, given by $F^{\sharp}(f):=
F\circ f$, being well-defined. So one sees that any
A$_\infty$-bimodule-map induces a map between the corresponding
Hochschild-cochain-complexes.

Finally, in section 5, for a given A$_\infty$-algebra $(A,D)$ one
can canonically make $A$ and its dual $A^{*}$ into an
A$_\infty$-bimodule over $A$ by using the map $D$. An
A$_\infty$-bimodule map between $A$ and $A^{*}$ will be called an
$\infty$-inner-product on $A$. There is a nice combinatorial
interpretation for $\infty $-inner-products, which will give rise
to a graph-complex associated to A$_\infty$-algebras with
$\infty$-inner-product.

The reason for the naming of the above definitions is that (just
like for the A$_\infty$-structure) one can write down lots of maps
which uniquely determine the A$_\infty$-bimodule-structure [etc.],
and which are a direct generalization of the usual concepts of
bimodule-structures [etc.].

All the spaces ($V$, $W$, $Z$, $A$, $M$, $N$, ...) in this paper
are always understood to be graded modules $V=\bigoplus_{i\in
\mathbb{Z}}V_{i}$ over a given ground ring $R$. The degree of
homogeneous elements $v\in V_{i}$ is written as $|v|:=i$, and the
degree of maps $\varphi:V_{i} \longrightarrow W_{j}$ is written as
$|\varphi|:=j-i$. All tensor-products of maps and their
compositions are understood in a graded way:
$$ (\varphi\otimes\psi)(v\otimes w)= (-1)^{|\psi|\cdot|v|}
   (\varphi(v))\otimes(\psi(w)), $$
$$ (\varphi\otimes\psi)\circ(\chi\otimes\varrho)= (-1)^{|\psi|\cdot|\chi|}
   (\varphi\circ\chi)\otimes(\psi\circ\varrho). $$
All objects $a_{i}$, $v_{i}$,... are assumed to be elements in
$A$, $V$,... respectively, if not stated otherwise (e.g.
Proposition 4.4.).\\
It will be necessary to look at elements of $V^{\otimes i}\otimes
V^{\otimes j}$. In order to distinguish between the tensor-product
in $V^{\otimes i}$ and the one between $V^{\otimes i}$ and
$V^{\otimes j}$, it is convenient to write the first one as a
tuple $(v_{1},...,v_{i})\in V^{\otimes i}$, and then
$(v_{1},...,v_{i})\otimes (v'_{1},...,v'_{j}) \in V^{\otimes
i}\otimes V^{\otimes j}$. The total degree of $(v_{1},...,v_{i})
\in V^{\otimes i}$ is given by $ |(v_{1},...,v_{i})|:=
\sum_{k=1}^{i} |v_{k}|$. \\
Frequently there will be sums of the form $ \sum_{i=0}^{n}
(v_{1},...,v_{i})\otimes(v_{i+1},...,v_{n}).$ Here the convention
will be used that for $i=0$, one has the term $1\otimes
(v_{1},...,v_{n})$ and for $i=n$ the term in the sum is
$(v_{1},...,v_{n})\otimes 1$, with $1=1_{TV}\in TV$. Similar for
terms $\sum_{i=0}^{k}(v_{1},...,v_{i}) \otimes(v_{i+1},...,v_{k},
w,v_{k+1},...,v_{n})$ the expression for $i=k$ is understood to
mean $(v_{1},...,v_{k}) \otimes (w,v_{k+1},...,v_{n})$.

I would like to thank Dennis Sullivan, who suggested to me that
there should be a sequence of homotopies describing
Poincar\'e-duality on the chain level, and M. Markl for many
useful comments.

\section{A$_\infty$-algebras}

Let us review the usual definitions about A$_\infty$-algebras that
are important to the discussion of this paper, and that can be
found in many sources (compare \cite{GJ} section 1, and
\cite{S2}).

\begin{defn}
A \textbf{coalgebra} $(C,\Delta)$ over a ring $R$ consists of an
$R$-module $C$ and a comultiplication  $\Delta:C\longrightarrow
C\otimes C$ of degree $0$ satisfying coassociativity:
\[
\begin{diagram}
\node{C}\arrow{e,t}{\Delta} \arrow{s,l}{\Delta}
\node{C\otimes C}\arrow{s,r}{\Delta\otimes id}\\
\node{C\otimes C}\arrow{e,b}{id\otimes \Delta} \node{C\otimes
C\otimes C}
\end{diagram}
\]
Then a \textbf{coderivation} on $C$ is a map $f:C\longrightarrow
C$ such that
\[
\begin{diagram}
\node{C}\arrow{e,t}{\Delta} \arrow{s,l}{f}
\node{C\otimes C}\arrow{s,r}{f\otimes id+id\otimes f}\\
\node{C}\arrow{e,b}{\Delta} \node{C\otimes C}
\end{diagram}
\]
\end{defn}

\begin{defn} Let $V=\bigoplus_{j\in \mathbb{Z}} V_{j}$ be a graded
module over a given ground ring $R$. The \textbf{tensor-coalgebra}
of $V$ over the ring $R$ is given by
$$ TV:=\bigoplus_{i\geq 0} V^{\otimes i}, $$
$$ \Delta:TV\longrightarrow TV\otimes TV, \,\,\,\,\,\,\,\,\,
   \Delta(v_{1},...,v_{n}):=\sum_{i=0}^{n} (v_{1},...,v_{i})
    \otimes(v_{i+1},...,v_{n}).$$
Let $A=\bigoplus_{j\in \mathbb{Z}} A_{j}$ be a graded module over
a given ground ring $R$. Define its \textbf{suspension} $sA$ to be
the graded module $sA= \bigoplus_{j\in \mathbb{Z}} (sA)_{j}$ with
$(sA)_{j}:= A_{j-1}$. The suspension map $s:A\longrightarrow sA$,
$s:a\mapsto
sa:=a$ is an isomorphism of degree +1.\\
Now the \textbf{bar complex} of $A$ is given by $BA:=T(sA)$.\\
An \textbf{A$_\infty$-algebra} on $A$ is given by a coderivation
$D$ on $BA$ of degree $-1$ such that $D^{2}=0$.
\end{defn}
Let's try to understand this definition.

The tensor-coalgebra has the property to lift every module map
$f:TV\longrightarrow V$ to a coalgebra-map $F:TV \longrightarrow
TV$:
\[
\begin{diagram}
\node{} \node{TV}\arrow{s,r}{projection}\\
\node{TV}\arrow{e,b}{f}\arrow{ne,t}{F} \node{V}
\end{diagram}
\]
A similar property for coderivations on $TV$ will make it possible
to understand the definition of A$_\infty$-algebras in a different
way.
\begin{lem}
\begin{itemize}
\item [(a)] Given a map $\varrho:V^{\otimes n}\longrightarrow V$
of degree $|\varrho|$, which can be viewed as a map
$\varrho:TV\longrightarrow V$ by letting its only non-zero
component being given by the original $\varrho$ on $V^{\otimes
n}$. Then $\varrho$ lifts uniquely to a coderivation
$\tilde{\varrho} :TV \longrightarrow TV$ with
\[
\begin{diagram}
\node{} \node{TV}\arrow{s,r}{projection} \\
\node{TV}\arrow{ne,t}{\tilde{\varrho}} \arrow{e,b}{\varrho}
\node{V}
\end{diagram}
\]
by taking
$$ \tilde{\varrho}(v_{1},...,v_{k}):=0, \,\,\,\,\,\,\, for \,\,\, k<n, $$
$$ \tilde{\varrho}(v_{1},...,v_{k}):=\sum_{i=0}^{k-n}
   (-1)^{|\varrho|\cdot(|v_{1}|+...+|v_{i}|)}(v_{1},...,\varrho
   (v_{i+1},...,v_{i+n}),...,v_{k}), $$
$$ \,\,\,\,\,\,\,\,\,\,\,\,\,\,\,\,\,\,\,\,\,\,\,\,\,\,\,\,\,\,\,\,
\,\,\,\,\,\,\,\,\,\,\,\,\,\,\,\,\,\,\,\,\,\,\,\,\,\,\,\,\,\,\,\,\,
\,\,\,\,\,\,\,\,\,\,\,\,\,\,\,\,\,\,\,\,\,\,\,\,\,\,\,\,\,\,\,\,\,
\,\,\,\,\,\,\,\,\,\,\,\,\,\,\,\,\,\,\,\,\,\,\,\,\,\,\,\,\,\,\,\,\,
\,\,\,\,\,\,\,\,\,\,\,\,\,\,\,\,\,\,\,\,\,\,\,\,\,\,\,\,\,\,\,\,\,
for \,\,\, k\geq n. $$ Thus $\tilde{\varrho}\mid_{V^{\otimes
k}}:V^{\otimes k} \longrightarrow V^{\otimes k-n+1} $.
\item [(b)] There is a one-to-one correspondence between
coderivations $\sigma:TV\longrightarrow TV$ and systems of maps
$\{\varrho_{i}:V^{\otimes i}\longrightarrow V\}_{i\geq 0}$, given
by $\sigma=\sum_{i\geq 0} \tilde{\varrho_{i}}$.
\end{itemize}
\end{lem}
\begin{proof}
\begin{itemize}
\item [(a)] The argument here is dual to the way one lifts
derivations on $TA$. To be precise one should use induction on the
output-component of $\tilde{\varrho}$. Denote by
$\tilde{\varrho}^{j}$ the component of $\tilde{\varrho}$ mapping
$TV\longrightarrow V^{\otimes j}$. Then
$\tilde{\varrho}^{1},...,\tilde{\varrho}^ {m-1}$ determine
uniquely the component $\tilde{\varrho}^{m}$, because of the
coderivation property of $\tilde{\varrho}$. Let's derive an
equation with which this can be seen easily.
\begin{eqnarray*}
\,\,\,\,\,\,\,\,\,\,\,\,\,\,
\Delta(\tilde{\varrho}(v_{1},...,v_{k}))&=&(\tilde{\varrho}\otimes
id+id\otimes \tilde{\varrho})(\Delta(v_{1},...,v_{k}))=\\
&=&(\tilde{\varrho}\otimes id+id\otimes \tilde{\varrho}) (
\sum_{i=0}^{k}(v_{1},...,v_{i})\otimes (v_{i+1},...,v_{k}))=\\
&=& \sum_{i=0}^{k}\tilde{\varrho}(v_{1},...,v_{i})\otimes
(v_{i+1},...,v_{k}) + \\
& & \,\,\, +(-1)^{|\tilde{\varrho}|\cdot(|v_{1}|+...+|v_{i}|)}
(v_{1},...,v_{i})\otimes \tilde{\varrho} (v_{i+1},...,v_{k}).
\end{eqnarray*}
Now, projecting both sides to $\bigoplus_{i+j=m}V^{\otimes i}
\otimes V^{\otimes j} \subset TV\otimes TV$ yields
\begin{eqnarray*}
\,\,\,\,\,\,\,\,\,\,\,\,\,\,
\Delta(\tilde{\varrho}^{m}(v_{1},...,v_{k}))&=&
\sum_{i=0}^{k}\tilde{\varrho}^{m+i-k}(v_{1},...,v_{i})\otimes
(v_{i+1},...,v_{k}) + \\
& & \,\,\, +(-1)^{|\tilde{\varrho}|\cdot(|v_{1}|+...+|v_{i}|)}
(v_{1},...,v_{i})\otimes \tilde{\varrho}^{m-i}
(v_{i+1},...,v_{k}).
\end{eqnarray*}
So the righthand side depends only on $\tilde{\varrho}^{j}$ with
$j<m$, except for the uninteresting terms $\tilde{\varrho}^{m}
(v_{1},...,v_{k})\otimes 1$ and $1 \otimes \tilde{\varrho}^{m}
(v_{1},...,v_{k})$, where $1\in TV$. With this, an induction
argument shows that $\tilde{\varrho}^{m}$ is only nonzero on
$V^{\otimes k}$ for $k=m+n-1$, where it is
$$ \tilde{\varrho}^{m}(v_{1},...,v_{m+n-1})=\sum_{i=0}^{m-1}
   (-1)^{|\varrho|\cdot(|v_{1}|+...+|v_{i}|)}(v_{1},...,\varrho
   (v_{i+1},...,v_{i+n}),...,v_{m+n-1}). $$
\item [(b)] Observe that being a coderivation is a linear
condition, and so the sum of coderivations is again a
coderivation. Thus the map
$$ \,\,\,\,\,\,\,\,\,\,\,\,\,\,\,\,\,\,\,\,\,\,\,
\alpha:\{\{\varrho_{i}:V^{\otimes i}\rightarrow V\}_{i\geq 0}\}
\longrightarrow Coder(TV), \,\,\,\,\,\,\, \{\varrho_{i}:
V^{\otimes i}\rightarrow V\}_{i\geq 0} \mapsto \sum_{i\geq 0}
\tilde{\varrho_{i}} $$ is well defined. Its inverse $\beta$ is
given by $\beta:\sigma\mapsto\{pr_{V}\circ\sigma|_{V^{\otimes i
}}\}_{i\geq 0}$, because the explicit lifting property of (a)
shows that $\beta\circ\alpha=id$, and the uniqueness part of (a)
shows that $\alpha\circ\beta=id$.
\end{itemize}
\end{proof}

Let's apply this to Definition 2.2.
\begin{prop} Let $(A,D)$ be an A$_\infty$-algebra. Now let $D$ be
given by a system of maps $\{D_{i}:sA^{\otimes i} \longrightarrow
sA\}_{i\geq 1}$, (where $D_{0}=0$ is assumed,) just like in Lemma
2.3.(b), and rewrite them as $m_{i}:A^{\otimes i} \longrightarrow
A$ given by $D_{i}=s\circ m_{i}\circ (s^{-1})^{\otimes i}$.\\
Then the condition $D^{2}=0$ is equivalent to the following system
of equations:
\begin{eqnarray*}
m_{1}(m_{1}(a_{1}))&=&0,\\
m_{1}(m_{2}(a_{1},a_{2}))-m_{2}(m_{1}(a_{1}),a_{2})-(-1)^{|a_{1}|}
m_{2}(a_{1},m_{1}(a_{2}))&=&0,\\
m_{1}(m_{3}(a_{1},a_{2},a_{3}))-m_{2}(m_{2}(a_{1},a_{2}),a_{3})+
m_{2}(a_{1},m_{2}(a_{2},a_{3}))+& &\\
+m_{3}(m_{1}(a_{1}),a_{2},a_{3})+(-1)^{|a_{1}|}m_{3}
(a_{1},m_{1}(a_{2}),a_{3})+& &\\
+(-1)^{|a_{1}|+|a_{2}|}m_{3}(a_{1},a_{2},m_{1}(a_{3}))&=&0,\\
...\\
\sum_{i=1}^{k} \sum_{j=0}^{k-i+1} (-1)^{\varepsilon} \cdot
m_{k-i+1} (a_{1},...,m_{i} (a_{j},...,a_{j+i-1}),...,a_{k})&=&0,\\
where \,\,\,\varepsilon=i\cdot \sum_{l=1}^{j-1}|a_{l}|+ (j-1)
\cdot(i+1)+k-i\\
...
\end{eqnarray*}
\end{prop}
\begin{proof}
The only difficulty is to determine the signs when replacing the
$D_{r}$'s by the $m_{r}$'s. First notice that
\begin{eqnarray*}
D_{k}(sa_{1},...,sa_{k})&=&s\circ m_{k}\circ (s^{-1})^{\otimes
k}(sa_{1},...,sa_{k})=\\
&=&(-1)^{\sum_{j=1}^{k-1}(|a_{j}|+1)} s\circ m_{k}\circ ((s^{-1})
^{\otimes k-1}\otimes id) (sa_{1},...,sa_{k-1},s^{-1}sa_{k})=\\
&=&(-1)^{\sum_{j=1}^{k-2}2(|a_{j}|+1)+|a_{k-1}|+1} s\circ
m_{k}\circ ((s^{-1})^{\otimes k-2}\otimes (id)^{\otimes 2})\\
& & \,\,\,\,\,\,\,\,\,\,\,\,\,\,\,\,\,\,\,\,\,\,\,\,\,\,\,\,\,\,
\,\,\,\,\,\,\,\,\,\,\,\,\,\,\,\,\,\,\,\,\,\,\,\,\,\,\,\,\,\,\,\,\,\,\,
 (sa_{1},...,sa_{k-2},s^{-1}sa_{k-1},s^{-1}sa_{k})=\\
& & \,\,\,\,\,\,\,\,\,\,\,\,\,\,\,\,\,\,\,\,\,\,\,\,\,\,\,\, ...\\
&=&(-1)^{\sum_{j=1}^{k}(k-j)\cdot(|a_{j}|+1)} s\circ m_{k}
(a_{1},...,a_{k}).
\end{eqnarray*}
Therefore it follows from Lemma 2.3.(a) applied to the degree $-1$
coderivation $D$, that
\begin{equation}
 pr_{sA}\circ D^{2}(sa_{1},...,sa_{k}) = \,\,\,\,\,\,\,\,\,\,\,\,
\,\,\,\,\,\,\,\,\,\,\,\,\,\,\,\,\,\,\,\,\,\,\,\,\,\,\,\,\,\,\,\,\,
\,\,\,\,\,\,\,\,\,\,\,\,\,\,\,\,\,\,\,\,\,\,\,\,\,\,\,\,\,\,\,\,\,
\,\,\,\,\,\,\,\,\,\,\,\,\,\,\,\,\,\,\,\,\,\,\,\,\,\,\,\,\,\,\,\,\,
\end{equation}
\begin{eqnarray*}
&=& \sum_{i=1}^{k} \sum_{j=0}^{k-i+1}(-1)^{\sum_{l=1}^{j-1}
(|a_{l}|+1)} D_{k-i+1} (sa_{1},...,D_{i}
(sa_{j},...,sa_{j+i-1}),...,sa_{k})=\\
&=& \sum_{i=1}^{k} \sum_{j=0}^{k-i+1}(-1)^{\sum_{l=1}^{j-1}
(|a_{l}|+1)+\sum_{l=j}^{j+i-1}(j+i-l-1)\cdot(|a_{l}|+1)} \cdot \\
& & \,\,\,\,\,\,\,\,\,\,\,\,\,\,\,\,\,\,\,\,\,\,\,\,\,\,\,
\,\,\,\,\,\,\,\,\,\,\,\,\,\,\,\,\,\,\,\,\,\,\,\,\,\,\,\,\,
\,\,\,\,\,\,\,\,\,\,\,\,\,\,\,\,\,\,\,
\cdot D_{k-i+1} (sa_{1},...,s\circ m_{i}(a_{j},...,a_{j+i-1}),...,sa_{k})=\\
&=& \sum_{i=1}^{k} \sum_{j=0}^{k-i+1}(-1)^{\varepsilon} m_{k-i+1}
(a_{1},...,m_{i} (a_{j},...,a_{j+i-1}),...,a_{k}),
\end{eqnarray*}
where $\varepsilon$ can be determined by using the fact that
$m_{r}=\pm s^{-1}\circ D_{r}\circ s^{\otimes r}$ is of degree
$-1-1+r=r-2$. Instead of doing the general case it is more
instructive to look at four special cases where $i$ and $j$ are
either even or odd. Let's take $k=8$. As seen above, signs occur
from applying $D_{i}$, transforming $D_{i}$ into $m_{i}$, and
transforming $D_{k-i+1}$ into $m_{k-i+1}$.\\
$i=3,\,\,\, j=3:$
\[
\begin{tabular} {c c c c c c c c c}
{ } & {$a_{1}$} & {$a_{2}$} & {($a_{3}$} & {$a_{4}$} & {$a_{5}$)} &
{$a_{6}$} & {$a_{7}$} & {$a_{8}$} \\ \hline
{ $D_{i}:$} & {$|a_{1}|+1$} & {$|a_{2}|+1$} & { } &
 { } & { } & { } & { } & { } \\ \hline
{ $m_{i}:$} & { } & { } & { } & {$|a_{4}|+1$} & { }
& { } & { } & { } \\ \hline
{ $m_{k-i+1}:$} & {$|a_{1}|+1$} & { } & {($|a_{3}|+$} &
{$|a_{4}|+|a_{5}|$} & {$+3-1$)} & { } & {$|a_{7}|+1$} & { } \\ \hline
\end{tabular}
\]
$i=3,\,\,\, j=4:$
\[
\begin{tabular} {c c c c c c c c c}
{ } & {$a_{1}$} & {$a_{2}$} & {$a_{3}$} & {($a_{4}$} & {$a_{5}$} &
{$a_{6}$)} & {$a_{7}$} & {$a_{8}$} \\ \hline
{ $D_{i}:$} & {$|a_{1}|+1$} & {$|a_{2}|+1$} & {$|a_{3}|+1$} &
 { } & { } & { } & { } & { } \\ \hline
{ $m_{i}:$} & { } & { } & { } & { } & {$|a_{5}|+1$}
& { } & { } & { } \\ \hline
{ $m_{k-i+1}:$} & {$|a_{1}|+1$} & { } & {$|a_{3}|+1$} & { } &
{ } & { } & {$|a_{7}|+1$} & { } \\ \hline
\end{tabular}
\]
So, if $i$ is odd, then the lower two rows ($m_{i}$ and
$m_{k-i+1}$) show that here the "$a_{r}$"-terms are exactly
$\sum_{l=1}^{k}(k-l)\cdot(|a_{l}|+1)$, and for $j=3$ there is an
additional $-1$. The top ($D_{i}$-)row has the terms
$\sum_{l=1}^{j-1}(|a_{l}|+1)=(\sum_{l=1}^{j-1}|a_{l}|)+j-1$. The
additional $-1$ can be put together with the $j-1$ to give a
constant depending only on $k$. Thus the term for $\varepsilon$ is
given by
$$ \varepsilon=\sum_{l=1}^{k}(k-l)\cdot(|a_{l}|+1)+ \sum_{l=1}^{j-1}|a_{l}|
+ k -1 $$ \\
$i=4,\,\,\, j=3:$
\[
\begin{tabular} {c c c c c c c c c}
{ } & {$a_{1}$} & {$a_{2}$} & {($a_{3}$} & {$a_{4}$} & {$a_{5}$} &
{$a_{6}$)} & {$a_{7}$} & {$a_{8}$} \\ \hline
{ $D_{i}:$} & {$|a_{1}|+1$} & {$|a_{2}|+1$} & { } &
 { } & { } & { } & { } & { } \\ \hline
{ $m_{i}:$} & { } & { } & {$|a_{3}|+1$} & { } & {$|a_{5}|+1$}
& { } & { } & { } \\ \hline
{ $m_{k-i+1}:$} & { } & {$|a_{2}|+1$} & { } & { } &
{ } & { } & {$|a_{7}|+1$} & { } \\ \hline
\end{tabular}
\]
$i=4,\,\,\, j=4:$
\[
\begin{tabular} {c c c c c c c c c}
{ } & {$a_{1}$} & {$a_{2}$} & {$a_{3}$} & {($a_{4}$} & {$a_{5}$} &
{$a_{6}$} & {$a_{7}$)} & {$a_{8}$} \\ \hline
{ $D_{i}:$} & {$|a_{1}|+1$} & {$|a_{2}|+1$} & {$|a_{3}|+1$} &
 { } & { } & { } & { } & { } \\ \hline
{ $m_{i}$:} & { } & { } & { } & {$|a_{4}|+1$} & { }
& {$|a_{6}|+1$} & { } & { } \\ \hline
{ $m_{k-i+1}$:} & { } & {$|a_{2}|+1$} & { } & {($|a_{4}|+$} &
{$|a_{5}|+|a_{6}|$} & {$+|a_{7}|+$} & {$4-1)$}  & { } \\ \hline
\end{tabular}
\]
So, if $i$ is even then the "$a_{r}$"-terms (from all three rows)
are exactly $\sum_{l=1}^{k}(k-l)\cdot(|a_{l}|+1)$. The only
$j$-dependence is the additional $-1$ in the $i=4,\,\,\, j=4$
case, which will induce an alternating sign (, which starts
depending on $k$). So, this case gives:
$$ \varepsilon=\sum_{l=1}^{k}(k-l)\cdot(|a_{l}|+1)+k+j-1. $$

Putting both cases together, and bringing the common term
$\sum_{l=1}^{k}(k-l)\cdot(|a_{l}|+1)$ (for both $i$ even and $i$
odd) to the left, we get the sign
\begin{eqnarray*}
 \varepsilon-(\sum_{l=1}^{k}(k-l)\cdot(|a_{l}|+1)) &=&
 i\cdot(\sum_{l=1}^{j-1}|a_{l}|+ k -1)+(i+1)\cdot(k+j-1)=\\
 &=& i\cdot\sum_{l=1}^{j-1}|a_{l}|+ (j-1)\cdot(i+1)+k-i.
\end{eqnarray*}
Thus, dividing the equation $D^{2}=0$ by the sign
$(-1)^{\sum_{l=1}^{k} (k-l)\cdot(|a_{l}|+1)}$ yields the result.
\end{proof}

\begin{expl} Any differential graded algebra $(A,\partial,\mu)$ gives an
A$_\infty$-algebra-structure on $A$ by taking $m_{1}:=\partial$,
$m_{2}:=\mu$ and $m_{k}:=0$ for $k\geq 3$. Then the equations from
Proposition 2.4. are the defining conditions of a differential
graded algebra:
\begin{eqnarray*}
\partial^{2}(a)&=&0,\\
\partial (a\cdot b)&=& \partial(a)\cdot b+ (-1)^{|a|}a\cdot \partial(b), \\
(a\cdot b)\cdot c &=& a\cdot (b\cdot c). \\
\end{eqnarray*}
There are no higher equations.
\end{expl}

\begin{defn} Given an A$_\infty$-algebra $(A,D)$. Then
the \textbf{Hochschild-cochain-complex of $A$} is defined to be
the space $C^{*}(A):=CoDer(BA,BA)$ of coderivations on $BA$ with
the differential $\delta:C^{*}(A)\longrightarrow C^{*}(A)$ given
by $\delta(f):=[D,f]=D\circ f-(-1)^{|f|}f\circ D$. It is
$\delta^{2}=0$, because as $D$ is of degree $-1$ and $D^{2}=0$ one
follows that $\delta^{2}(f)=[D,D\circ f-(-1)^{|f|}f\circ D]=
D\circ D\circ f -(-1)^{|f|} D\circ f\circ D -(-1)^{|f|+1} D\circ f
\circ D + (-1)^{|f|+|f|+1} f\circ D \circ D=0$.
\end{defn}

\section{A$_\infty$-bimodules}

Given an A$_\infty$-algebra $(A,D)$. The goal is now to define the
concept of an A$_\infty$-bimodule over $A$, which was already
given in \cite{GJ} and \cite{M}.\\
This should be a generalization of two facts. First, one can
define the  Hochschild-cochain-complex for any algebra with values
in a bimodule, and so one should still be able to make that
definition in the infinite case. Second, any algebra is a bimodule
over itself by left- and right-multiplication, which should again
still hold (see section 5).

The following space and map are important ingredients.

\begin{defn} For modules $V$ and $W$ over $R$, one writes
$$T^{W}V:=R\oplus\bigoplus_{k\geq 0, l\geq 0} V^{\otimes k} \otimes W
\otimes V^{\otimes l}.$$
Furthermore, let
$$ \Delta^{W}:T^{W}V\longrightarrow (TV\otimes T^{W}V)\oplus
(T^{W}V\otimes TV), $$ be given by
\begin{eqnarray*}
\Delta^{W}(v_{1},...,v_{k},w,v_{k+1},...,v_{k+l})&:=&
 \sum_{i=0}^{k}(v_{1},...,v_{i}) \otimes(v_{i+1},...,w,...,v_{n})+\\
 & &+\sum_{i=k}^{k+l} (v_{1},...,w,...,v_{i})
 \otimes(v_{i+1},...,v_{k+l}).
\end{eqnarray*}
Again for modules $A$ and $M$ let $B^{M}A$ be given by $T^{sM}sA$,
where $s$ is the suspension from Definition 2.2.
\end{defn}

Observe that $T^{W}V$ is not a coalgebra, but rather a bi-comodule
over $TV$. Here is a definition for a coderivation from $TV$ to
$T^{W}V$.

\begin{defn} A \textbf{coderivation} from $TV$ to $T^{W}V$ is defined to be a map
$f:TV\longrightarrow T^{W}V$ that makes the following diagram
commute:
\[
\begin{diagram}
\node{TV}\arrow{s,l}{f}\arrow{e,t}{\Delta}
  \node{TV\otimes TV}\arrow{s,r}{id\otimes f+f\otimes id} \\
\node{T^{W}V}\arrow{e,b}{\Delta^{W}} \node{(TV\otimes
T^{W}V)\oplus (T^{W}V\otimes TV)}
\end{diagram}
\]
For modules $A$ and $M$ let $C^{*}(A,M):=CoDer(BA,B^{M}A)$ be the
space of coderivations in the above sense. This space is called
the \textbf{Hochschild-cochain-complex of $A$ with values in $M$}.
\end{defn}

\begin{lem}
\begin{itemize}
\item [(a)]
Given a map $\varrho:V^{\otimes n}\longrightarrow W$ of degree
$|\varrho|$, which can be viewed as a map
$\varrho:TV\longrightarrow W$ by letting its only non-zero
component being given by the original $\varrho$ on $V^{\otimes
n}$. Then $\varrho$ lifts uniquely to a coderivation
$\tilde{\varrho} :TV \longrightarrow T^{W}V$ with
\[
\begin{diagram}
\node{} \node{T^{W}V}\arrow{s,r}{projection} \\
\node{TV}\arrow{ne,t}{\tilde{\varrho}} \arrow{e,b}{\varrho}
\node{W}
\end{diagram}
\]
by taking
$$ \tilde{\varrho}(v_{1},...,v_{k}):=0, \,\,\,\,\,\,\, for \,\,\, k<n, $$
$$ \tilde{\varrho}(v_{1},...,v_{k}):=\sum_{i=0}^{k-n}
   (-1)^{|\varrho|\cdot(|v_{1}|+...+|v_{i}|)}
   (v_{1},...,\varrho(v_{i+1},...,v_{i+n}),...,v_{k}),$$
$$ \,\,\,\,\,\,\,\,\,\,\,\,\,\,\,\,\,\,\,\,\,\,\,\,\,\,\,\,\,\,\,\,
\,\,\,\,\,\,\,\,\,\,\,\,\,\,\,\,\,\,\,\,\,\,\,\,\,\,\,\,\,\,\,\,\,
\,\,\,\,\,\,\,\,\,\,\,\,\,\,\,\,\,\,\,\,\,\,\,\,\,\,\,\,\,\,\,\,\,
\,\,\,\,\,\,\,\,\,\,\,\,\,\,\,\,\,\,\,\,\,\,\,\,\,\,\,\,\,\,\,\,\,
\,\,\,\,\,\,\,\,\,\,\,\,\,\,\,\,\,\,\,\,\,\,\,\,\,\,\,\,\,\,\,\,\,
for \,\,\, k\geq n. $$ Thus $\tilde{\varrho}\mid_{V^{\otimes
k}}:V^{\otimes k} \longrightarrow \bigoplus_{i+j=k-n} V^{\otimes
i} \otimes W \otimes V^{\otimes j}$.
\item [(b)] There is a one-to-one correspondence between
coderivations $\sigma:TV\longrightarrow T^{W}V$ and systems of
maps $\{\varrho_{i}:V^{\otimes i}\longrightarrow W\}_{i\geq 0}$,
given by $\sigma=\sum_{i\geq 0} \tilde{\varrho_{i}}$.
\end{itemize}
\end{lem}
\begin{proof}
\begin{itemize}
\item [(a)]
The proof is similar to the one of Lemma 2.3.(a). Now
$\tilde{\varrho}^{j}$ is meant to be the component of
$\tilde{\varrho}$ mapping $TV\longrightarrow\bigoplus_{r+s=j}
V^{\otimes r} \otimes W\otimes V^{\otimes s}$. Let's do again
induction on $m$ for $\tilde{\varrho}^{m}$. The equation
\begin{eqnarray*}
\,\,\,\,\,\,\,\,\,\,\,\,\,\,
\Delta^{W}(\tilde{\varrho}(v_{1},...,v_{k}))&=&(\tilde{\varrho}\otimes
id+id\otimes \tilde{\varrho})(\Delta(v_{1},...,v_{k}))=\\
&=&(\tilde{\varrho}\otimes id+id\otimes \tilde{\varrho}) (
\sum_{i=0}^{k}(v_{1},...,v_{i})\otimes (v_{i+1},...,v_{k}))=\\
&=& \sum_{i=0}^{k}\tilde{\varrho}(v_{1},...,v_{i})\otimes
(v_{i+1},...,v_{k}) + \\
& & \,\,\, +(-1)^{|\tilde{\varrho}|\cdot(|v_{1}|+...+|v_{i}|)}
(v_{1},...,v_{i})\otimes \tilde{\varrho} (v_{i+1},...,v_{k})
\end{eqnarray*}
is being projected to the component
$$ \bigoplus_{r+s+t=m}(V^{r}\otimes W\otimes V^{s})\otimes V^{t}+
V^{r}\otimes(V^{s}\otimes W\otimes V^{t}) \subset T^{W}V\otimes
TV+TV\otimes T^{W}V. $$ This gives
\begin{eqnarray*}
\,\,\,\,\,\,\,\,\,\,\,\,\,\,
\Delta^{W}(\tilde{\varrho}^{m}(v_{1},...,v_{k}))&=&
\sum_{i=0}^{k}\tilde{\varrho}^{m+i-k}(v_{1},...,v_{i})\otimes
(v_{i+1},...,v_{k}) + \\
& & \,\,\, +(-1)^{|\tilde{\varrho}|\cdot(|v_{1}|+...+|v_{i}|)}
(v_{1},...,v_{i})\otimes \tilde{\varrho}^{m-i}
(v_{i+1},...,v_{k}).
\end{eqnarray*}
The righthand side depends only on $\tilde{\varrho}^{j}$ with
$j<m$, except for the uninteresting terms $\tilde{\varrho}^{m}
(v_{1},...,v_{k})\otimes 1$ and $1 \otimes \tilde{\varrho}^{m}
(v_{1},...,v_{k})$. An induction shows that $\tilde{\varrho}^{m}$
is only nonzero on $V^{\otimes k}$ for $k=m+n-1$, where it is
$$ \tilde{\varrho}^{m}(v_{1},...,v_{m+n-1})=\sum_{i=0}^{m-1}
   (-1)^{|\varrho|\cdot(|v_{1}|+...+|v_{i}|)}(v_{1},...,\varrho
   (v_{i+1},...,v_{i+n}),...,v_{m+n-1}). $$
\item [(b)]
Then maps
\begin{eqnarray*}
& \alpha:\{\{\varrho_{i}:V^{\otimes i}\rightarrow W\}_{i\geq 0}\}
\longrightarrow Coder(TV,T^{W}V), & \{\varrho_{i}: V^{\otimes
i}\rightarrow W\}_{i\geq 0} \mapsto
\sum_{i\geq 0} \tilde{\varrho_{i}} \\
& \beta:Coder(TV,T^{W}V) \longrightarrow \{\{\varrho_{i}:
V^{\otimes i}\rightarrow W\}_{i\geq 0}\}, &
\sigma\mapsto\{pr_{W}\circ\sigma|_{V^{\otimes i }}\}_{i\geq 0}
\end{eqnarray*}
are inverse to each other by (a).
\end{itemize}
\end{proof}

Of course one wants to put a differential on $C^{*}(A,M)$ just
like in section 2.

\begin{prop} Given an A$_{\infty}$-algebra $(A,D)$ and a module $M$.
Let $D^{M}:B^{M}A\longrightarrow B^{M}A$ be a map of degree $-1$.
Then the induced map $\delta^{M}:CoDer(BA,B^{M}A) \longrightarrow
CoDer(BA,B^{M}A)$, given by $\delta^{M}(f):=D^{M}\circ
f-(-1)^{|f|}f\circ D$, is well-defined, (i.e. it maps
coderivations to coderivations,) if and only if the following
diagram commutes:
\begin{equation}
\begin{diagram}
\node{B^{M}A}\arrow{s,l}{D^{M}}\arrow{e,t}{\Delta^{M}}
  \node{(BA\otimes B^{M}A)\oplus (B^{M}A\otimes BA)}\arrow{s,r}
  {(id\otimes D^{M}+D\otimes id) \oplus (D^{M}\otimes id+id\otimes D)} \\
\node{B^{M}A}\arrow{e,b}{\Delta^{M}} \node{(BA\otimes
  B^{M}A)\oplus (B^{M}A\otimes BA)}
\end{diagram}
\end{equation}
\end{prop}
\begin{proof} Let $f:BA\longrightarrow B^{M}A$ be a coderivation.
One needs to investigate under which conditions $\delta^{M}(f)$ is
also a coderivation. This means, that
$$ (id\otimes\delta^{M}(f)+\delta^{M}(f)\otimes id)\circ\Delta=
\Delta^{M}\circ\delta^{M}(f), $$ or
$$ (id\otimes(D^{M}\circ f)-(-1)^{|f|}id\otimes(f\circ D)+
    (D^{M}\circ f)\otimes id-(-1)^{|f|}
    (f\circ D)\otimes id)\circ\Delta=\,\,\,\, $$
$$ \,\,\,\,\,\,\,\,\,\,\,\,\,\,\,\,\,\,\,\,\,\,\,\,\,\,\,\,\,
   \,\,\,\,\,\,\,\,\,\,\,\,\,\,\,\,\,\,\,\,\,\,\,\,\,\,\,\,\,
    \,\,\,\,\,\,\,\,\,= \Delta^{M}\circ D^{M}\circ f-(-1)^{|f|}
    \Delta^{M}\circ f \circ D. $$
Now, using the coderivation property for $f$ and $D$, one gets the
following identity
\begin{eqnarray*}
\Delta^{M}\circ f \circ D &=&
 (id\otimes f)\circ\Delta\circ D + (f\otimes id)\circ\Delta\circ D=\\
 &=& (id\otimes f)\circ(id\otimes D)\circ \Delta+
     (id\otimes f)\circ(D\otimes id)\circ \Delta+\\
 &&  +(f\otimes id)\circ(id\otimes D)\circ \Delta+
     (f\otimes id)\circ(D\otimes id)\circ \Delta=\\
 &=& (id\otimes(f\circ D) +(-1)^{|f|} D\otimes f + f\otimes D +
      (f\circ D)\otimes id) \circ \Delta,
\end{eqnarray*}
and so the requirement for $\delta^{M}(f)$ above reduces to
\begin{eqnarray*}
\Delta^{M}\circ D^{M}\circ f&=& (id\otimes(D^{M}\circ f)+
(D^{M}\circ f)\otimes id+ D\otimes f +(-1)^{|f|} f\otimes D)\circ \Delta=\\
&=& (id\otimes D^{M}+ D\otimes id)\circ(id\otimes f)\circ \Delta +
    (D^{M}\otimes id+ id\otimes D)\circ(f\otimes id)\circ \Delta=\\
&=& (id\otimes D^{M}+ D\otimes id)\circ\Delta^{M} \circ f +
    (D^{M}\otimes id+ id\otimes D)\circ\Delta^{M} \circ f.\\
\end{eqnarray*}
The last step looks strange, because $\Delta^{M} \circ f=
(id\otimes f)\circ \Delta+ (f\otimes id)\circ \Delta$. But as
$id\otimes D^{M}+ D\otimes id:BA\otimes B^{M}A \longrightarrow
BA\otimes B^{M}A$, this map doesn't pick up any part from
$(f\otimes id)\circ \Delta:BA \longrightarrow B^{M}A \otimes BA$.
A similar argument applies to $D^{M}\otimes id+ id
\otimes D$.\\
So, $D^{M}$ has to satisfy
$$ \Delta^{M}\circ D^{M}\circ f = (id\otimes D^{M}+ D\otimes id+
    D^{M}\otimes id+ id\otimes D)\circ\Delta^{M} \circ f, $$
for all coderivations $f:TA\longrightarrow T^{M}A$. By Lemma 3.3.,
one sees that there are enough coderivations, to make this
condition equivalent to
$$ \Delta^{M}\circ D^{M}= (id\otimes D^{M}+ D\otimes id+
    D^{M}\otimes id+ id\otimes D)\circ\Delta^{M},$$
which is the claim.
\end{proof}

Again one wants to describe $D^{M}$ by a system of maps.

\begin{lem}
\begin{itemize}
\item [(a)] Given a module $V$ and a coderivation $\psi$ on $TV$
 with an associated system of maps $\{\psi_{i}:V^{\otimes i} \longrightarrow
 V\}_{i\geq 1}$ from Lemma 2.3., where $\psi_{i}$ is of degree
 $|\psi_{i}|$. Then any map $\varrho:T^{W}V\longrightarrow W$
 given by $\varrho=\sum_{k\geq0, l\geq0} \varrho_{k,l}$, with
 $\varrho_{k,l}: V^{\otimes k}\otimes W \otimes V^{\otimes l}
 \longrightarrow W$ of degree $|\varrho_{k,l}|$, lifts uniquely to a map
 $\tilde{\varrho} :T^{W}V \longrightarrow T^{W}V$
\[
\begin{diagram}
\node{} \node{T^{W}V}\arrow{s,r}{projection} \\
\node{T^{W}V}\arrow{e,b}{\varrho} \arrow{ne,t}{\tilde{\varrho}}
\node{W}
\end{diagram}
\]
which makes the following diagram commute (compare diagram (3.1)
in Proposition 3.4.)
\begin{equation}
\begin{diagram}
\node{T^{W}V}\arrow{s,l}{\tilde{\varrho}}\arrow{e,t}{\Delta^{W}}
  \node{(TV\otimes T^{W}V)\oplus (T^{W}V\otimes TV)}\arrow{s,r}
  {(id\otimes \tilde{\varrho}+\psi\otimes id) \oplus (\tilde{\varrho}
  \otimes id+id\otimes \psi)} \\
\node{T^{W}V}\arrow{e,b}{\Delta^{W}} \node{(TV\otimes
  T^{W}V)\oplus (T^{W}V\otimes TV)}
\end{diagram}
\end{equation}
This map is given by taking
$$ \tilde{\varrho}(v_{1},...,v_{k},w,v_{k+1},...,v_{k+l}) :=
\,\,\,\,\,\,\,\,\,\,\,\,\,\,\,\,\,\,\,\,\,\,\,\,\,\,\,\,\,\,
\,\,\,\,\,\,\,\,\,\,\,\,\,\,\,\,\,\,\,\,\,\,\,\,\,\,\,\,\,\,
\,\,\,\,\,\,\,\,\,\,\,\,\,\,\,\,\,\,\,\,\,\,\,\,\,\,\,\,\,\,$$
\begin{eqnarray*}
&:=& \sum_{i=1}^{k} \sum_{j=1}^{k-i+1}
     (-1)^{|\psi_{i}|\sum_{r=1}^{j-1}|v_{r}|}
     (v_{1},...,\psi_{i}(v_{j},...,v_{i+j-1}),...,w,...,v_{k+l}) +\\
&&   +\sum_{i=0}^{k} \sum_{j=0}^{l}
     (-1)^{|\varrho_{i,j}|\sum_{r=1}^{k-i}|v_{r}|}
     (v_{1},...,\varrho_{i,j}(v_{k-i+1},...,w,...,v_{k+j}),...,v_{k+l}) +\\
&&   +\sum_{i=1}^{l} \sum_{j=1}^{l-i+1}
     (-1)^{|\psi_{i}|(|w|+\sum_{r=1}^{k+j-1}|v_{r}|)}
     (v_{1},...,w,...,\psi_{i}(v_{k+j},...,v_{k+i+j-1}),...,v_{k+l}).
\end{eqnarray*}
(Notice that the condition of satisfying diagram (3.2) is not
linear, i.e. if $\chi$ and $\psi$ both make diagram (3.2) commute,
then $\chi+\psi$ will not.)
\item [(b)]
There is a one-to-one correspondence between maps
$\sigma:T^{W}V\longrightarrow T^{W}V$ that make diagram (3.2)
commute and maps $\varrho=\sum\varrho_{k,l}$ like in (a), given by
$\sigma=\tilde{\varrho}$.
\end{itemize}
\end{lem}
\begin{proof}
\begin{itemize}
\item [(a)]
Again one uses induction on the output-component of
$\tilde{\varrho}$. Denote by $\tilde{\varrho}^{j}$ the component
of $\tilde{\varrho}$ mapping $T^{W}V\longrightarrow \bigoplus
_{k+l=j} V^{\otimes k}\otimes W\otimes V^{\otimes l}$ and by
$\psi^{j}$ the component of $\psi$ mapping $TV\longrightarrow
V^{\otimes}$. Then $\tilde{\varrho}^{0},...,\tilde{\varrho}^
{m-1}$ will uniquely determine the component $\tilde
{\varrho}^{m}$.
$$ \Delta^{W}(\tilde{\varrho}(v_{1},...,v_{k},w,v_{k+1},...,
   v_{k+l}))= \,\,\,\,\,\,\,\,\,\,\,\,\,\,\,\,\,\,\,\,\,
   \,\,\,\,\,\,\,\,\,\,\,\,\,\,\,\,\,\,\,\,\,\,\,\,\,\,\, $$
\begin{eqnarray*}
&=&(id\otimes \tilde{\varrho}+\psi\otimes id+ \tilde{\varrho}
   \otimes id+id\otimes\psi) (\Delta^{W}
   (v_{1},...,v_{k},w,v_{k+1},...,v_{k+l}))=\\
&=&(id\otimes \tilde{\varrho}+\psi\otimes id+ \tilde{\varrho}
   \otimes id+id\otimes\psi)
   (\sum_{i=0}^{k}(v_{1},...,v_{i}) \otimes (v_{i+1},...,w,...,v_{k+l})+ \\
& & \,\,\,\,\,\,\,\,\,\,\,\,\,\,\,\,\,\,\,\,\,\,\,\,\,\,\,\,\,\,
\,\,\,\,\,\,\,\,\,\,\,\,\,\,\,\,\,\,\,\,\,\,\,\,\,\,\,\,\,\,\,\,
\,\,\,\,\,\,\,\,\,\,\,\,\,\,\,\,\,\,\,\,\,\,\,\,\,\,\,\,\,\,\,\,
    +\sum_{i=k}^{k+l}(v_{1},...,w,...,v_{i}) \otimes (v_{i+1},...,v_{k+l}))=\\
&=& \sum_{i=0}^{k}(-1)^{|\tilde{\varrho}|\sum_{r=1}^{i}|v_{r}|}
    (v_{1},...,v_{i}) \otimes \tilde{\varrho}(v_{i+1},...,w,...,v_{k+l})+\\
& & +\sum_{i=0}^{k}\psi(v_{1},...,v_{i}) \otimes
(v_{i+1},...,w,...,v_{k+l})+
    \sum_{i=k}^{k+l}\tilde{\varrho}(v_{1},...,w,...,v_{i})
    \otimes (v_{i+1},...,v_{k+l})+\\
& & +\sum_{i=k}^{k+l} (-1)^{|\psi|(|w|+\sum_{r=1}^{i}|v_{r}|)}
    (v_{1},...,w,...,v_{i}) \otimes \psi(v_{i+1},...,v_{k+l}). \\
\end{eqnarray*}
Then projecting both sides to
$$ \bigoplus_{r+s+t=m}(V^{\otimes r} \otimes W\otimes V^{\otimes s})
\otimes V^{\otimes t}+ V^{\otimes } \otimes (V^{\otimes s} \otimes
W\otimes V^{\otimes t}) \subset T^{W}V\otimes TV + TV\otimes
T^{W}V$$ yields
\begin{eqnarray*}
\,\,\,\,\,\,\,\,\,\,\,\,\,\,
\Delta^{W}(\tilde{\varrho}^{m}(v_{1},...,w,...,v_{k}))&=&
    \sum_{i=0}^{k}\pm (v_{1},...,v_{i}) \otimes
    \tilde{\varrho}^{m-i}(v_{i+1},...,w,...,v_{k+l})+\\
& & +\sum_{i=0}^{k}\psi^{m+i-k-l}(v_{1},...,v_{i}) \otimes
    (v_{i+1},...,w,...,v_{k+l})+\\
& & +\sum_{i=k}^{k+l}\tilde{\varrho}^{m+i-k-l}
    (v_{1},...,w,...,v_{i}) \otimes (v_{i+1},...,v_{k+l})+\\
& & +\sum_{i=k}^{k+l} \pm (v_{1},...,w,...,v_{i})
    \otimes \psi^{m-i}(v_{i+1},...,v_{k+l}). \\
\end{eqnarray*}
So the righthand side depends only on $\psi^{j}$'s, which are all
explicitly known by Lemma 2.3., and $\tilde{\varrho}^{j}$ with
$j<m$, (except for the uninteresting terms $\tilde{\varrho}^{m}
(v_{1},...,w,... ,v_{k+l})\otimes 1$ and $1 \otimes
\tilde{\varrho}^{m} (v_{1},...,w,...,v_{k+l})$). With this, an
induction argument shows that $\tilde{\varrho}^{m}$ is given by
the formula of the Lemma.
\item [(b)]
Let $X:=\{\sigma:T^{W}V\longrightarrow T^{W}V\,\,|\,\, \sigma$
makes diagram (3.2) commute $\}$. Then
\begin{eqnarray*}
\alpha:\{\varrho:T^{W}V\longrightarrow W\}\longrightarrow X,
& & \varrho \mapsto \tilde{\varrho}, \\
\beta:X\longrightarrow \{\varrho:T^{W}V\longrightarrow W\}, & &
\sigma \mapsto pr_{W}\circ\sigma
\end{eqnarray*}
are inverse to each other by (a).
\end{itemize}
\end{proof}

\begin{defn} Given an A$_{\infty}$-algebra $(A,D)$. Then an
\textbf{A$_{\infty}$-bimodule} $(M,D^{M})$ consists of a module
$M$ together with a map $D^{M}:B^{M}A\longrightarrow B^{M}A$ of
degree $-1$, which makes the diagram (3.1) of Proposition 3.4.
commute,
and satisfies $(D^{M})^{2} =0$.\\
By Proposition 3.4. one can put the differential $\delta^{M}:
CoDer(TA,T^{M}A) \longrightarrow CoDer(TA,T^{M}A)$, $\delta(f)
:=D^{M}\circ f-(-1)^{|f|}f\circ D$ on the
Hochschild-cochain-complex. Now it satisfies $(\delta^{M})^{2}=0$,
because with $(D^{M})^{2} =0$, one gets
$(\delta^{M})^{2}(f)=D^{M}\circ D^{M}\circ f-(-1)^{|f|}D^{M}\circ
f\circ D -(-1)^{|f|+1}
D^{M}\circ f \circ D + (-1)^{|f|+|f|+1} f\circ D\circ D =0$.\\
The definition of an A$_{\infty}$-bimodule was already given in
\cite{GJ} section 3 and also in \cite{M}, and coincides with the
one here.
\end{defn}

\begin{prop} Let $(A,D)$ be an A$_\infty$-algebra with a system of
maps $\{m_{i}:A^{\otimes i} \longrightarrow A\}_{i\geq 1}$
associated to $D$ by Proposition 2.4. (where $m_{0}=0$ is
assumed). Let $(M,D^{M})$ be an A$_\infty$-bimodule over $A$ with
a system of maps $\{D^{M}_{k,l} :sA^{\otimes k}\otimes sM \otimes
sA^{\otimes l} \longrightarrow sM\}_{k\geq 0, l\geq 0}$ from Lemma
3.5.(b) associated to $D^{M}$. Let $b_{k,l}:A^{\otimes k}\otimes M
\otimes A^{\otimes l} \longrightarrow M$ be the induced map by
$D^{M}_{k,l}=s\circ b_{k,l}\circ (s^{-1})^{\otimes k+l+1}$.\\
Then the condition $(D^{M})^{2}=0$ is equivalent to the following
system of equations:
\begin{eqnarray*}
b_{0,0}(b_{0,0}(m))&=&0, \\
b_{0,0}(b_{0,1}(m,a_{1}))-b_{0,1}(b_{0,0}(m),a_{1})-(-1)^{|m|}
b_{0,1}(m,m_{1}(a_{1})) &=& 0,\\
b_{0,0}(b_{1,0}(a_{1},m))-b_{1,0}(m_{1}(a_{1}),m)-(-1)^{|a_{1}|}
b_{1,0}(a_{1},b_{0,0}(m)) &=& 0,\\
b_{0,0}(b_{1,1}(a_{1},m,a_{2}))-b_{0,1}(b_{1,0}(a_{1},m),a_{2})+
b_{1,0}(a_{1},b_{0,1}(m,a_{2}))+& &\\
+b_{1,1}(m_{1}(a_{1}),m,a_{2})
+(-1)^{|a_{1}|} b_{1,1}(a_{1},b_{0,0}(m),a_{2})+& &\\
+(-1)^{|a_{1}|+|m|}b_{1,1}(a_{1},m,m_{1}(a_{2})) &=& 0,\\
\end{eqnarray*}
$$...$$
\begin{eqnarray*}
\sum_{i=1}^{k} \sum_{j=1}^{k-i+1} \pm
 b_{k-i+1,l}(a_{1},...,m_{i}(a_{j},...,a_{i+j-1}),...,m,...,a_{k+l}) +& &\\
+\sum_{i=0}^{k} \sum_{j=0}^{l} \pm
 b_{k-i,l-j}(a_{1},...,b_{i,j}(a_{k-i+1},...,m,...,a_{k+j}),...,a_{k+l}) +& &\\
+\sum_{i=1}^{l} \sum_{j=1}^{l-i+1} \pm
 b_{k,l-i+1}(a_{1},...,m,...,m_{i}(a_{k+j},...,a_{k+i+j-1}),...,a_{k+l})&=&0
\end{eqnarray*}
$$...$$
\end{prop}
where the signs are exactly analogous to the ones in Proposition
2.4.
\begin{proof}
The result follows immediately from Lemma 3.5., after rewriting
$D^{M}_{k,l}$ and $D_{j}$ by $b_{k,l}$ and $m_{j}$. (Notice that
this replacement only changes a sign.)\\
In order to get the correct sign, first notice that the lifting
described in Lemma 3.5. is exactly the usual lifting as
coderivations, except that one has to pick $b_{k,l}$ or $m_{j}$
according to where the element of $M$ is located. Keeping this in
mind, it is possible to redo all the steps from Proposition 2.4.
\end{proof}

\begin{expl} Let's pick up Example 2.5. Let $(A,\partial,\mu)$ be a
differential graded algebra with the A$_\infty$-algebra-structure
$m_{1}:=\partial$, $m_{2}:=\mu$ and $m_{k}:=0$ for $k\geq 3$. Now,
let $(M,\partial',\lambda,\rho)$ be a differential graded bimodule
over $A$, where $\lambda:A\otimes M\longrightarrow M$ and $\rho
:M\otimes A\longrightarrow M$ denote the left- and right-action.
It is possible to make $M$ into an A$_\infty$-bimodule over $A$ by
taking $b_{0,0}:=\partial'$, $b_{1,0}:=\lambda$, $b_{0,1}:=\rho$
and $b_{k,l}:=0$ for $k+l>1$. Then the equations of Proposition
3.7. are the defining conditions of a differential bialgebra over
$A$:
\begin{eqnarray*}
(\partial')^{2}(m)&=&0,\\
\partial'(m.a)&=&m.\partial(a)+(-1)^{|m|}\partial'(m).a,\\
\partial'(a.m)&=&\partial(a).m+(-1)^{|a|}a.\partial'(m),\\
(a.m).b&=&a.(m.b),\\
(m.a).b&=&m.(a\cdot b),\\
a.(b.m)&=&(a\cdot b).m.
\end{eqnarray*}
There are no higher equations.
\end{expl}

For later purposes it is convenient to have the following

\begin{lem} Given an A$_\infty$-algebra $(A,D)$ and an
A$_\infty$-bimodule $(M,D^{M})$, with system of maps $\{b_{k,l}
:A^{\otimes k}\otimes M \otimes A^{\otimes l} \longrightarrow
M\}_{k\geq 0, l\geq 0}$ from Proposition 3.7.\\
Then the dual space $M^{*}:=Hom_{R}(M,R)$ has a canonical
A$_\infty$-bimodule-structure given by maps $\{b'_{k,l}
:A^{\otimes k}\otimes M^{*} \otimes A^{\otimes l} \longrightarrow
M^{*}\}_{k\geq 0, l\geq 0}$,
$$ (b'_{k,l}(a_{1},...,a_{k},m^{*},a_{k+1},...,a_{k+l}))(m):=(-1)^
   {\varepsilon}\cdot m^{*}(b_{l,k}(a_{k+1},...,a_{k+l},m,a_{1},...,a_{k})), $$
$$ \,\,\,\,\,\,\,\,\,\,\,\,\,\,\,\,\,\,\,\,
where \,\,\,\varepsilon:= (|a_{1}|+...+|a_{k}|)\cdot
(|m^{*}|+|a_{k+1}|+...+|a_{k+l}|+|m|)+|m^{*}|\cdot(k+l+1)
.$$
\end{lem}
\begin{proof} First, it is easy to see that the degrees work out,
because of $|b'_{k,l}|=|b_{l,k}|$.\\
In order to show $(D^{M^{*}})^{2}=0$, one can use the criterion
from Proposition 3.7. The top and the bottom term in the general
sum of Proposition 3.7. convert to
$$ (b'_{k-i+1,l}(a_{1},...,m_{i}(a_{j},...,a_{i+j-1})
    ,...,m^{*},...,a_{k+l}))(m)=
\,\,\,\,\,\,\,\,\,\,\,\,\,\,\,\,\,\,\,\,\,\,\,\,\,\,\,\,\,\,\,\,\,$$
$$ \,\,\,\,\,\,\,\,\,\,  =\pm m^{*}(b_{l,k-i+1}(a_{k+1},...,a_{k+l},m,
a_{1},...,m_{i}(a_{j},...,a_{i+j-1}),...,a_{k})), $$ and
$$ (b'_{k,l-i+1}(a_{1},...,m^{*},...,m_{i}(a_{k+j},...,
   a_{k+i+j-1}),...,a_{k+l}))(m)=
\,\,\,\,\,\,\,\,\,\,\,\,\,\,\,\,\,\,\,\,\,\,\,\,\,\,\,\,\,\,\,\,\,$$
$$ \,\,\,\,\,\,\,\,\,\,  =\pm m^{*}(b_{l-i+1,k}(a_{k+1},...,m_{i}(a_{k+j},...,
   a_{k+i+j-1}),...,a_{k+l},m,a_{1},...,a_{k})). $$
So, these terms come from terms of the
A$_\infty$-bimodule-structure of $M$. The same is true for the
middle term:
$$ (b'_{k-i,l-j}(a_{1},...,b'_{i,j}(a_{k-i+1},...,
 m^{*},...,a_{k+j}),...,a_{k+l}))(m) =
\,\,\,\,\,\,\,\,\,\,\,\,\,\,\,\,\,\,\,\,\,\,\,\,\,\,\,\,\,\,\,\,\,$$
\begin{eqnarray*}
&=& \pm (b'_{i,j}(a_{k-i+1},...,m^{*},...,a_{k+j}))(b_{l-j,k-i}
      (a_{k+j+1},...,a_{k+l},m,a_{1},...,a_{k-i}))=\\
&=& \pm m^{*}(b_{j,i}(a_{k+1},...,a_{k+j},b_{l-j,k-i}
      (a_{k+j+1},...,a_{k+l},m,a_{1},...,a_{k-i}),a_{k-i+1},...,a_{k})).
\end{eqnarray*}
So, the sum from Proposition 3.7. for the A$_\infty$-bimodule
$M^{*}$ has exactly the terms of $m^{*}$ applied the the sum for
the A$_\infty$-bimodule $M$.\\
The only remaining question is whether the signs are correct. The
proof for this is left to the reader.
\end{proof}

\section{Morphisms of A$_\infty$-bimodules}

Given two A$_\infty$-bimodules $(M,D^{M})$ and $(N,D^{N})$ over an
A$_\infty$-algebra $(A,D)$. What is the natural notion of morphism
between them?

Again a motivation is to have for any A$_\infty$-bimodule-map an
induced map of their Hochschild-cochain-complexes.

\begin{prop} Given three modules $V$, $W$ and $Z$. Let $F:T^{W}V
\longrightarrow T^{Z}V$ be a map. Then the induced map $F^{\sharp}
:CoDer(TV,T^{W}V) \longrightarrow CoDer(TV,T^{Z}V)$, given by
$F^{\sharp}(f):=F\circ f$, is well-defined, (i.e. it maps
coderivations to coderivations,) if and only if the following
diagram commutes:
\begin{equation}
\begin{diagram}
\node{T^{W}V}\arrow{s,l}{F}\arrow{e,t}{\Delta^{W}}
  \node{(TV\otimes T^{W}V)\oplus (T^{W}V\otimes TV)}\arrow{s,r}
  {(id\otimes F) \oplus (F\otimes id)} \\
\node{T^{Z}V}\arrow{e,b}{\Delta^{Z}} \node{(TV\otimes
  T^{Z}V)\oplus (T^{Z}V\otimes TV)}
\end{diagram}
\end{equation}
\end{prop}
\begin{proof} If both $f:TV\longrightarrow T^{W}V$ and $F\circ
f:TV\longrightarrow T^{Z}V$ are coderivations then this means that
the top diagram and the overall diagram below commute.
\[
\begin{diagram}
\node{TV}\arrow{s,l}{f}\arrow{e,t}{\Delta}
  \node{TV\otimes TV}\arrow{s,r} {(id\otimes f) + (f\otimes id)} \\
\node{T^{W}V}\arrow{s,l}{F}\arrow{e,t}{\Delta^{W}}
  \node{(TV\otimes T^{W}V)\oplus (T^{W}V\otimes TV)}\arrow{s,r}
  {(id\otimes F) \oplus (F\otimes id)} \\
\node{T^{Z}V}\arrow{e,b}{\Delta^{Z}} \node{(TV\otimes
  T^{Z}V)\oplus (T^{Z}V\otimes TV)}
\end{diagram}
\]
But then the lower diagram has to commute if applied to any
element in $Im(f)\subset T^{W}V$. By Lemma 3.3. there are enough
coderivations to make this true for all $T^{W}V$.
\end{proof}

Again let's describe $F$ by a system of maps.

\begin{lem}
\begin{itemize}
\item [(a)] Given modules $V$, $W$ and $Z$ and a map
 $\varrho:V^{\otimes k}\otimes W\otimes V^{\otimes l}\longrightarrow Z$
 of degree $|\varrho|$, which can be viewed as a map $\varrho:
 T^{W}V \longrightarrow Z$ by letting its only nonzero component
 be  the original $\varrho$ on $V^{\otimes k}\otimes W\otimes V^{\otimes
 l}$. Then $\varrho$ lifts uniquely to a map $\tilde{\varrho} :T^{W}V
 \longrightarrow T^{Z}V$
\[
\begin{diagram}
\node{} \node{T^{Z}V}\arrow{s,r}{projection} \\
\node{T^{W}V} \arrow{e,b}{\varrho} \arrow{ne,t}{\tilde{\varrho}}
\node{Z}
\end{diagram}
\]
which makes the diagram (4.1) in Proposition 4.1. commute (put
$\tilde{\varrho}$ instead of $F$). This map is given by
$$ \tilde{\varrho}(v_{1},...,v_{r},w,v_{r+1},...,v_{r+s}):=0,
   \,\,\,\,\,\,\, for \,\,\, r<k \,\,\, or \,\,\, s<l,$$
$$ \tilde{\varrho}(v_{1},...,v_{r},w,v_{r+1},...,v_{r+s}):=
 \,\,\,\,\,\,\,\,\,\,\,\,\,\,\,\,\,\,\,\,\,\,\,\,\,\,\,\,\,\,
 \,\,\,\,\,\,\,\,\,\,\,\,\,\,\,\,\,\,\,\,\,\,\,\,\,\,\,\,\,\,$$
$$  \,\,\,\,\,\,\,\,\,\,\,\,\,\,\,\,\,\,\,\,\,\,\,\,\,\,\,\,\,\,
   := (-1)^{|\varrho|\sum_{i=1}^{r-k}|v_{i}|}
   (v_{1},...,\varrho(v_{r-k+1},...,w,...,v_{r+l}),...,v_{r+s}),$$
$$ \,\,\,\,\,\,\,\,\,\,\,\,\,\,\,\,\,\,\,\,\,\,\,\,\,\,\,\,\,\,\,
 \,\,\,\,\,\,\,\,\,\,\,\,\,\,\,\,\,\,\,\,\,\,\,\,\,\,\,\,\,\,\,
 \,\,\,\,\,\,\,\,\,\,\,\,\,\,\,\,\,\,\,\,\,\,\,\,\,\,\,\,\,\,\,
 \,\,\,\,\,\,\,\,\,\,\,\,\,\,\,\,\,\,\,\,\,\,\,\,\,\,\,\,\,\,\,
 \,\,\,\,\,\,\,\,\,\,\,\,\,\,\,\,\,\,\,\,\,\,\,\,\,\,\,\,\,\,\,
 for \,\,\, r\geq k \,\,\, and \,\,\, s\geq l. $$
Thus $\tilde{\varrho}\mid_{V^{\otimes r}\otimes W\otimes
V^{\otimes s}}:V^{\otimes r}\otimes W\otimes V^{\otimes s}
\longrightarrow V^{\otimes r-k} \otimes Z\otimes V^{\otimes s-l}$.
\item [(b)] There is a one-to-one correspondence between
maps $\sigma:T^{W}V\longrightarrow T^{Z}V$ making diagram (4.1)
commute and systems of maps $\{\varrho_{k,l}:V^{\otimes k}\otimes
W\otimes V^{\otimes l}\longrightarrow Z\}_{k\geq 0, l\geq 0}$,
given by $\sigma=\sum_{k\geq 0, l\geq 0}
\widetilde{\varrho_{k,l}}$.
\end{itemize}
\end{lem}
\begin{proof}
\begin{itemize}
\item [(a)]
Again one uses induction on the output-component of
$\tilde{\varrho}$. Denote by $\tilde{\varrho}^{j}$ the component
of $\tilde{\varrho}$ mapping $T^{W}V\longrightarrow \bigoplus
_{r+s=j} V^{\otimes r}\otimes Z\otimes V^{\otimes s}$. Then
$\tilde{\varrho}^{1},...,\tilde{\varrho}^ {m-1}$ determine
uniquely the component $\tilde{\varrho}^{m}$.
$$ \Delta^{Z}(\tilde{\varrho}(v_{1},...,v_{r},w,v_{r+1},...,
   v_{r+s}))= \,\,\,\,\,\,\,\,\,\,\,\,\,\,\,\,\,\,\,\,\,
\,\,\,\,\,\,\,\,\,\,\,\,\,\,\,\,\,\,\,\,\,\,\,\,\,\,\,\,\,\,\,\,
   \,\,\,\,\,\,\,\,\,\,\,\,\,\,\,\,\,\,\,\,\,\,\,\,\,\,\, $$
\begin{eqnarray*}
&=&(id\otimes \tilde{\varrho}+\tilde{\varrho} \otimes id)
   (\Delta^{W} (v_{1},...,v_{r},w,v_{r+1},...,v_{r+s}))=\\
&=&(id\otimes \tilde{\varrho}+\tilde{\varrho} \otimes id)
   (\sum_{i=0}^{r}(v_{1},...,v_{i}) \otimes (v_{i+1},...,w,...,v_{r+s})+ \\
& & \,\,\,\,\,\,\,\,\,\,\,\,\,\,\,\,\,\,\,\,\,\,\,\,\,\,\,\,\,\,
\,\,\,\,\,\,\,\,\,\,\,\,\,\,\,\,\,\,\,\,\,\,\,\,\,\,\,\,\,\,\,\,
    +\sum_{i=r}^{r+s}(v_{1},...,w,...,v_{i}) \otimes (v_{i+1},...,v_{r+s}))=\\
&=& \sum_{i=0}^{r}(-1)^{|\tilde{\varrho}|\sum_{t=1}^{i}|v_{t}|}
    (v_{1},...,v_{i}) \otimes \tilde{\varrho}(v_{i+1},...,w,...,v_{r+s})+\\
& & \,\,\,\,\,\,\,\,\,\,\,\,\,\,\,\,\,\,\,\,\,\,\,\,\,\,\,\,\,\,
\,\,\,\,\,\,\,\,\,\,\,\,\,\,\,\,\,\,\,\,\,\,\,\,\,\,\,\,\,\,\,\,
    +\sum_{i=r}^{r+s}\tilde{\varrho}(v_{1},...,w,...,v_{i})
    \otimes (v_{i+1},...,v_{r+s}).\\
\end{eqnarray*}
Then projecting both sides to
$$ \bigoplus_{r+s+t=m}(V^{\otimes r} \otimes Z\otimes V^{\otimes s})
\otimes V^{\otimes t}+ V^{\otimes t} \otimes (V^{\otimes s}
\otimes Z\otimes V^{\otimes t}) \subset T^{Z}V\otimes TV +
TV\otimes T^{Z}V $$ yields
\begin{eqnarray*}
\,\,\,\,\,\,\,\,\,\,\,\,\,\,
\Delta^{Z}(\tilde{\varrho}^{m}(v_{1},...,w,...,v_{r+s}))&=&
    \sum_{i=0}^{r}\pm (v_{1},...,v_{i}) \otimes
    \tilde{\varrho}^{m-i}(v_{i+1},...,w,...,v_{r+s})+\\
& & +\sum_{i=r}^{r+s}\tilde{\varrho}^{m+i-r-s}
    (v_{1},...,w,...,v_{i}) \otimes (v_{i+1},...,v_{r+s}).\\
\end{eqnarray*}
So the righthand side depends only on $\tilde{\varrho}^{j}$ with
$j<m$, (except for the uninteresting terms $\tilde{\varrho}^{m}
(v_{1},...,w,... ,v_{r+s})\otimes 1$ and $1 \otimes
\tilde{\varrho}^{m} (v_{1},...,w,...,v_{r+s})$). With this, an
induction argument shows that $\tilde{\varrho}^{m}$ is only
nonzero on $V^{\otimes r}\otimes W\otimes V^{\otimes s}$ with
$r-k+s-l=m$, where it is given by
$$ \tilde{\varrho}^{m}(v_{1},...,v_{r},w,v_{r+1},...,v_{r+s})=
 \,\,\,\,\,\,\,\,\,\,\,\,\,\,\,\,\,\,\,\,\,\,\,\,\,\,\,\,\,\,
 \,\,\,\,\,\,\,\,\,\,\,\,\,\,\,\,\,\,\,\,\,\,\,\,\,\,\,\,\,\,$$
$$  \,\,\,\,\,\,\,\,\,\,\,\,\,\,\,\,\,\,\,\,\,\,\,\,\,\,\,\,\,\,
   = (-1)^{|\varrho|\sum_{i=1}^{r-k}|v_{i}|}
   (v_{1},...,\varrho(v_{r-k+1},...,w,...,v_{r+l}),...,v_{r+s}). $$
\item [(b)]
Let $X:=\{\sigma:T^{W}V\longrightarrow T^{Z}V\,\,|\,\, \sigma$
makes diagram (4.1) commute $\}$. Then
\begin{eqnarray*}
& & \alpha:\{\{\varrho_{k,l}:V^{\otimes k}\otimes W\otimes
V^{\otimes l}\longrightarrow Z\}_{k\geq 0, l\geq 0}\}
\longrightarrow X, \\
& & \,\,\,\,\,\,\,\,\,\,\,\,\,\,\,\,\,\,\,
\{\varrho_{k,l}:V^{\otimes k}\otimes W\otimes V^{\otimes
l}\longrightarrow Z\}_{k\geq 0, l\geq 0} \mapsto
\sum_{k\geq 0, l\geq 0} \widetilde{\varrho_{k,l}},\\
& & \beta:X \longrightarrow \{\{\varrho_{k,l}:V^{\otimes k}\otimes
W\otimes V^{\otimes l}\longrightarrow Z\}_{k\geq 0, l\geq 0}\},\\
& & \,\,\,\,\,\,\,\,\,\,\,\,\,\,\,\,\,\,\,
\sigma\mapsto\{pr_{Z}\circ\sigma|_{V^{\otimes k}\otimes W\otimes
V^{\otimes l}}\}_{k\geq 0, l\geq 0}
\end{eqnarray*}
are inverse to each other by (a).
\end{itemize}
\end{proof}

Let's apply this to the Hochschild-space.

\begin{defn} Given two A$_\infty$-bimodules $(M,D^{M})$ and
$(N,D^{N})$ over an A$_\infty$-algebra $(A,D)$. Then a map $F:
B^{M}A\longrightarrow B^{N}A$ of degree 0 is called an
\textbf{A$_\infty$-bimodule-map} $:\Leftrightarrow$ $F$ makes the
diagram
\[
\begin{diagram}
\node{B^{M}A}\arrow{s,l}{F}\arrow{e,t}{\Delta^{M}}
  \node{(BA\otimes B^{M}A)\oplus (B^{M}A\otimes BA)}\arrow{s,r}
  {(id\otimes F) \oplus (F\otimes id)} \\
\node{B^{N}A}\arrow{e,b}{\Delta^{N}} \node{(BA\otimes
  B^{N}A)\oplus (B^{N}A\otimes BA)}
\end{diagram}
\]
commute, and in addition it holds that $F\circ D^{M}=D^{N}\circ F$.\\
By Proposition 4.1., every A$_\infty$-bimodule-map induces (by
composition $F^{\sharp}:f\mapsto F\circ f$) a map between the
Hochschild-spaces, which preserves the differentials, because
$(F^{\sharp}\circ\delta^{M})(f)=F^{\sharp}(D^{M}\circ f+(-1)^{|f|}
f\circ D)= F\circ D^{M}\circ f+(-1)^{|f|}F\circ f\circ D= D^{N}
\circ F\circ f+(-1)^{|f|} F\circ f\circ D=\delta^{N}(F\circ
f)=(\delta^{N}\circ F^{\sharp})(f)$.
\end{defn}

\begin{prop} Let $(A,D)$ be an A$_\infty$-algebra with a system of
maps $\{m_{i}:A^{\otimes i} \longrightarrow A\}_{i\geq 1}$ from
Proposition 2.4. associated to $D$, (where $m_{0}=0$ is assumed).
Let $(M,D^{M})$ and $(N,D^{N})$ be A$_\infty$-bimodules over $A$
with systems of maps $\{b_{k,l} :A^{\otimes k}\otimes M \otimes
A^{\otimes l} \longrightarrow M\}_{k\geq 0, l\geq 0}$ and
$\{c_{k,l} :A^{\otimes k}\otimes N \otimes A^{\otimes l}
\longrightarrow N\}_{k\geq 0, l\geq 0}$ from Proposition 3.7.
associated to $D^{M}$ and $D^{N}$ respectively. Let $F:T^{M}A
\longrightarrow T^{N}A$ be an A$_\infty$-bimodule-map between $M$
and $N$, and let $\{F_{k,l} :sA^{\otimes k}\otimes sM \otimes
sA^{\otimes l} \longrightarrow sN\}_{k\geq 0, l\geq 0}$ be a
system of maps associated to $F$ by Lemma 4.2.(b). Again, rewrite
the maps $F_{k,l}$ by $f_{k,l} :A^{\otimes k}\otimes M \otimes
A^{\otimes l} \longrightarrow N$ by using the suspension map:
$F_{k,l}=s\circ f_{k,l}\circ (s^{-1})^{\otimes k+l+1}$.\\
Then the condition $F\circ D^{M}=D^{N}\circ F$ is equivalent to
the following system of equations:
$$ f_{0,0}(b_{0,0}(m))=c_{0,0} (f_{0,0}(m)), $$
$$f_{0,0}(b_{0,1}(m,a))- f_{0,1}(b_{0,0}(m),a)- (-1)^{|m|}
f_{0,1}(m,m_{1}(a))= \,\,\,\,\,\,\,\,\,\,\,\,\,\,\,\,\,\,\,\,\,\,
\,\,\,\,\,\,\,\,\,\,\,\,\,\,\,\,\,\,\,\,\,\,$$
$$ \,\,\,\,\,\,\,\,\,\,\,\,\,\,\,\,\,\,\,\,\,\,\,\,\,\,\,\,\,\,\,\,\,
\,\,\,\,\,\,\,\,\,\,\,\,\,\,\,\,\,\,\,\,\,\,\,\,\,\,\,\,\,\,\,\,\,
\,\,\,\,\,\,\,\,\,=c_{0,0}(f_{0,1}(m,a))+ c_{0,1}(f_{0,0}(m),a),$$
$$f_{0,0}(b_{1,0}(a,m))- f_{1,0}(m_{1}(a),m)- (-1)^{|a|}
f_{1,0}(a,b_{0,0}(m))= \,\,\,\,\,\,\,\,\,\,\,\,\,\,\,\,\,\,\,\,\,\
\,\,\,\,\,\,\,\,\,\,\,\,\,\,\,\,\,\,\,\,\,\,$$
$$ \,\,\,\,\,\,\,\,\,\,\,\,\,\,\,\,\,\,\,\,\,\,\,\,\,\,\,\,\,\,\,\,\,
\,\,\,\,\,\,\,\,\,\,\,\,\,\,\,\,\,\,\,\,\,\,\,\,\,\,\,\,\,\,\,\,\,
\,\,\,\,\,\,\,\,\,=c_{0,0}(f_{1,0}(a,m))+ c_{1,0}(a,f_{0,0}(m)),$$
$$...$$
\begin{eqnarray*}
&\sum_{i=1}^{k} \sum_{j=1}^{k-i+1} (-1)^{\varepsilon}
 f_{k-i+1,l}(a_{1},...,m_{i}(a_{j},...,a_{i+j-1}),...,m,...,a_{k+l+1}) +& \\
&+\sum_{j=1}^{k} \sum_{i=k-j+2}^{k+l-j+2} (-1)^{\varepsilon}
 f_{j,k+l-i-j+3} \,\,\,\,\,\,\,\,\,\,\,\,\,\,\,\,\,\,\,\,
 \,\,\,\,\,\,\,\,\,\,\,\,\,\,\,\,\,\,\,\,\,\,\,\,\,\,\,\,\,\,
 \,\,\,\,\,\,\,\,\,\,\,\,\,\,\,\,\,\,\,\,\,\,\,\,\,\,\,\,\,\, &\\
& \,\,\,\,\,\,\,\,\,\,\,\,\,\,\,\,\,\,\,\,\,\,\,\,\,\,\,\,\,\,
(a_{1},...,b_{k-j+1,i+j-k-2}(a_{j},...,m,...,a_{i+j-1}),...,a_{k+l+1}) +& \\
&+\sum_{i=1}^{l} \sum_{j=k+2}^{k+l-i+2} (-1)^{\varepsilon}
 f_{k,l-i+1}(a_{1},...,m,...,m_{i}(a_{j},...,a_{i+j-1}),...,a_{k+l+1})&= \\
=&\sum_{j=1}^{k+1} \sum_{i=k-j+2}^{k+l-j+2} (-1)^{\varepsilon'}
 c_{j,k+l-i-j+3} \,\,\,\,\,\,\,\,\,\,\,\,\,\,\,\,\,\,\,\,
 \,\,\,\,\,\,\,\,\,\,\,\,\,\,\,\,\,\,\,\,\,\,\,\,\,\,\,\,\,\,
 \,\,\,\,\,\,\,\,\,\,\,\,\,\,\,\,\,\,\,\,\,\,\,\,\,\,\,\,\,\, &\\
& \,\,\,\,\,\,\,\,\,\,\,\,\,\,\,\,\,\,\,\,\,\,\,\,\,\,\,\,\,\,
(a_{1},...,f_{k-j+1,i+j-k-2}(a_{j},...,m,...,a_{i+j-1}),...,a_{k+l+1}) & \\
\end{eqnarray*}
In order to simplify notation, it is assume that in $(a_{1},...
,a_{k+l+1})$ above, only the first $k$ and the last $l$ elements
are elements of $A$ and $a_{k+1}=m\in M$. Then the signs are given
by
$$ \,\,\,\varepsilon=i\cdot \sum_{r=1}^{j-1}|a_{r}|+ (j-1)
\cdot(i+1)+(k+l+1)-i, $$
$$ and \,\,\,\varepsilon'=(i+1)\cdot (j+1+\sum_{r=1}^{j-1}|a_{r}|).
\,\,\,\,\,\,\,\,\,\,\,\,\,\,\,\,\,\,\,\,\,\,\,\,\,
\,\,\,\,\,\,\,\,\,\,\, $$
$$...$$
\end{prop}
\begin{proof}
Up to signs these formulas follow immediately from the explicit
lifting properties in Lemma 3.5.(a) and Lemma 4.2.(a). For the
sign, the arguments of Proposition 2.4. will be applied.\\
Let's assume again that in $(a_{1},... ,a_{k+l+1})$, only the
first $k$ and the last $l$ elements are elements of $A$ and
$a_{k+1}\in M$. Now notice that just like in Proposition 2.4. one
gets
$$ F_{k,l}(sa_{0},...,sa_{k+l})=(-1)^{\sum_{j=1}^{k+l+1}(k+l+1-j)
\cdot(|a_{j}|+1)} s\circ f_{k,l} (a_{0},...,a_{k+l}).$$

So, when writing out the term $pr_{sN}\circ F\circ D^{M}
(sa_{0},... ,sa_{k+l})$, exactly the same signs appear, that were
in equation (2.1) of Proposition 2.4. for $pr_{sA}\circ D\circ D$.
This is so, because $D^{M}$, which has to be applied in the
argument of $F$, has degree $-1$ just like $D$, and the
application of the suspension map is the same for $D$ or $D^{M}$
or $F$. It follows that in this case the signs can simply be taken
from Proposition 2.4.

Unfortunately the signs for the term $D^{N}\circ F$ cannot be
taken directly from Proposition 2.4. like above. The difference is
that $F$, which is of degree $0$ (and not $-1$), has to be applied
in the argument of $D^{N}$. So, when $F$ "jumps" over elements
$sa_{i}$, no signs are introduced. This means that here one gets a
difference in signs compared to Proposition 2.4. given by
$$ (-1)^{\sum_{r=1}^{j-1}(|a_{r}+1|)}=
   (-1)^{\sum_{r=1}^{j-1}|a_{r}|+j-1} $$
(compare this with the first equality in (2.1)). Here one has to
take the same interpretation for the variables $i$ and $j$ as in
Proposition 2.4.; namely $f_{r,s}$ takes exactly $i$ inputs and
the first variable in $f_{r,s}$ is given by $a_{j}$:
$$ (a_{1},...,f_{r,s}(a_{j},...,m,...,a_{i+j-1}),...,a_{k+l+1}). $$
Another difference to Proposition 2.4. is given in the last step
of equation (2.1), because the interior element $m_{i}(a_{j},...,
a_{j+i-1})$ is replaced by some $f_{r,s}(a_{j},..., a_{j+i-1})$,
($r+s=i-1$), with $|m_{i}|=i-2$ and $|f_{r,s}|=i-1$. So, when
converting $D_{k-i+1}$ to $m_{k-i+1}$ in (2.1), the suspension map
for $a_{j+i},...,a_{k}$ jumps over one degree less. In the given
case, this introduces a difference in signs of $(-1)^{k+l+1-i-j+1}$.\\
Putting this together with the sign in Proposition 2.4. gives
$$ \varepsilon-(\sum_{r=1}^{k+l+1}(k+l+1-r)\cdot(|a_{r}|+1)) =
\,\,\,\,\,\,\,\,\,\,\,\,\,\,\,\,\,\,\,\,\,\,\,\,\,\,\,\,\,\,\,
\,\,\,\,\,\,\,\,\,\,\,\,\,\,\,\,\,\,\,\,\,\,\,\,\,\,\,\,\,\,\,
\,\,\,\,\,\,\,\,\,\,\,\,\,\,\,\,\,\,\,\,\,\,\,\,\,\,\,\,\,\,\, $$
\begin{eqnarray*}
 &=& i\cdot\sum_{r=1}^{j-1}|a_{r}|+ (j-1)\cdot(i+1)+(k+l+1)-i+\\
 & & +(\sum_{r=1}^{j-1}|a_{r}|+j-1)-(k+l+1-i-j+1) \equiv \\
 &\equiv& (i+1)\cdot(j+1+\sum_{r=1}^{j-1}|a_{r}|)
 \,\,\,\,\,\,\,\,\,\,\,\,\,\,\,\,\,\,\,\,\, (mod\,\,\, 2) \\
\end{eqnarray*}
Thus, dividing the equation $D^{N}\circ F=F\circ D^{M}$ by the
sign $(-1)^{\sum_{r=1}^{k+l+1} (k+l+1-r)\cdot(|a_{r}|+1)}$ yields
the result.
\end{proof}

\begin{expl}
Let's pick up the examples 2.5. and 3.8. Let $(A,\partial, \mu)$
be a differential graded algebra with the
A$_\infty$-algebra-structure $m_{1}:=\partial$, $m_{2}:=\mu$ and
$m_{k}:=0$ for $k\geq 3$. Now, let $(M,\partial^{M},\lambda^{M},
\rho^{M})$ and $(N,\partial^{N},\lambda^{N},\rho^{N})$ be
differential graded bimodules over $A$, with the A$_\infty
$-bialgebra-structures given by $b_{0,0}:=\partial^{M}$,
$b_{1,0}:=\lambda^{M}$, $b_{0,1}:=\rho^{M}$ and $b_{k,l}:=0$ for
$k+l>1$, and $c_{0,0}:=\partial^{N}$, $c_{1,0}:=\lambda^{N}$,
$c_{0,1}:=\rho^{N}$ and $c_{k,l}:=0$ for $k+l>1$.\\
Given a bialgebra map $f:M\longrightarrow N$ of degree 0. Then one
makes $f$ into a map of A$_\infty $-bialgebras by taking
$f_{0,0}:=f$ and $ f_{k,l}:=0$ for $k+l>0$. Then the equations of
Proposition 4.4. are the defining conditions of a differential
bialgebra map from $M$ to $N$:
\begin{eqnarray*}
f\circ\partial^{M} (m) &=& \partial^{N}\circ f(m)\\
f(m.a) &=& f(m).a\\
f(a.m) &=& a.f(m)
\end{eqnarray*}
There are no higher equations.
\end{expl}

\section{$\infty$-inner-products on A$_\infty$-algebras}

There are canonical A$_\infty$-bialgebra-structures on a given
A$_\infty$-algebra and its dual. A$_\infty$-bialgebra-maps between
them will then be defined to be $\infty$-inner products.

\begin{lem} Given an A$_\infty$-algebra $(A,D)$. Let the
coderivation $D$ be given by the system of maps $\{m_{i}:
A^{\otimes i} \longrightarrow A\}_{i\geq 1}$ from Proposition 2.4.
\begin{itemize}
\item [(a)] One can define an A$_\infty$-bimodule-structure on $A$
by taking $b_{k,l}:A^{\otimes k}\otimes A \otimes A^{\otimes
l}\longrightarrow A$ to be given by
$$ b_{k,l}:=m_{k+l+1}.$$
\item [(b)] One can define an A$_\infty$-bimodule-structure on $A^{*}$
by taking $b_{k,l}:A^{\otimes k}\otimes A^{*} \otimes A^{\otimes
l}\longrightarrow A^{*}$ to be given by
$$ (b_{k,l}(a_{1},...,a_{k},a^{*},a_{k+1},...,a_{k+l}))(a):=
\,\,\,\,\,\,\,\,\,\,\,\,\,\,\,\,\,\,\,\,\,\,\,\,\,\,\,\,\,\,\,\,\,\,\,$$
$$\,\,\,\,\,\,\,\,\,\,\,\,\,\,\,\,\,\,\,\,\,\,\,\,\,\,\,\,\,\,\,\,\,\,\,
   :=\pm a^{*}(m_{k+l+1}(a_{k+1},...,a_{k+l},a,a_{1},...,a_{k})),$$
where the signs are given in Lemma 3.9.
\end{itemize}
\end{lem}
\begin{proof}
\begin{itemize}
\item [(a)]  First notice that the A$_\infty$-bialgebra extension
described in Lemma 3.5.(a) becomes in this case the same as the
extension by coderivation described in Lemma 2.3.(a). Now, the
equations of Proposition 3.7. become exactly those of Proposition
2.4. and the diagram (3.1) from Proposition 3.4. becomes the usual
coderivation diagram for $D$.
\item [(b)] This follows immediately from (a) and Lemma 3.9.
\end{itemize}
\end{proof}

\begin{expl}
In the case of a differential algebra $(A,\partial, \mu)$, which
by Example 2.5. can be seen as an A$_\infty$-algebra, the above
A$_\infty$-bialgebra structure on $A$ is exactly the bialgebra
structure given by left- and right-multiplication, because then
$b_{1,0}(a\otimes b)=m_{2}(a\otimes b)=a\cdot b$ and
$b_{0,1}(a\otimes b)=m_{2}(a\otimes b)=a\cdot b$, for $a,b\in A$.\\
Similar the A$_\infty$-bialgebra structure on $A^{*}$ is given by
right- and left-multiplication in the arguments: $b_{1,0}(a\otimes
b^{*})(c)=b^{*}(m_{2}(c\otimes a))=b^{*}(c\cdot a)$ and
$b_{0,1}(a^{*}\otimes b)(c)=a^{*}(m_{2}(b\otimes c))=a^{*}(b\cdot
c)$, for $a,b,c\in A$, and $a^{*},b^{*}\in A^{*}$.
\end{expl}

\begin{defn} Given an A$_\infty$-algebra $(A,D)$. Then define an
\textbf{$\infty$-inner-product} on $A$ to be an
A$_\infty$-bimodule-map from the A$_\infty$-bimodule $A$ to the
A$_\infty$-bimodule $A^{*}$ given in Lemma 5.1.
\end{defn}

\begin{prop} Given an A$_\infty$-algebra $(A,D)$. Then an
$\infty$-inner product on $A$ is exactly given by a system of
inner-products on $A$, namely $\{<.,.,...>_{k,l}:A^{\otimes k+l+2}
\longrightarrow R\}_{k\geq 0, l\geq 0}$, that satisfies the
following relations:
$$ \sum_{i=1}^{k+l+2}(-1)^{\sum_{j=1}^{i-1}|a_{j}|}
 <a_{1},...,\partial(a_{i}),...,a_{k+l+2}>_{k,l}=\sum_{i,j,n}
 \pm <a_{i},...,m_{j}(a_{n},...),...>_{r,s}, $$
where in the sum on the right side, there is exactly
\underline{one} multiplication $m_{j}$ ($j\geq 2$) inside the
inner-product $<...>_{r,s}$ and this sum is taken over all i, j, n
subject to the following conditions:
\begin{itemize}
\item [(i)] The cyclic order of the $(a_{1},...,a_{k+l+2})$ is
preserved.
\item [(ii)] $a_{k+l+2}$ is always in the last slot of
$<...>_{r,s}$.
\item [(iii)] It might happen that $a_{k+l+2}$ is inside
$m_{j}$. By ii), this is the only case, when the inner product can
start with an $a_{i}\neq a_{1}$, \\ (e.g. $<a_{i+1},...,m_{j}
(a_{n},...,a_{k+l+2},a_{1},...,a_{i})>_{r,s}$ for $i\geq 1$).
\item [(iv)] $a_{k+1}$ and $a_{k+l+2}$ are never inside the $m_{j}$
together. (This is exactly the significance of the indices $k$ and
$l$.)
\item [(v)] $r$ and $s$ are given by looking at which slot the element
$a_{k+1}$ ends up in the inner-product. More exactly, $a_{k+1}$
will sit in the $(r+1)$st spot of $<...>_{r,s}$. $s$ is then
determined by saying that $<...>_{r,s}$ takes exactly $r+s+2$
arguments.
\end{itemize}
\end{prop}
\begin{proof} Let's use the description given in Proposition 4.4.
for A$_\infty$-bimodule-maps. An A$_\infty$-bimodule-map from $A$
to $A^{*}$ is given by maps $f_{k,l}:A^{\otimes k}\otimes A
\otimes A^{\otimes l}\longrightarrow A^{*}$, for $k,\,\,\, l\geq
0$. These can clearly be interpreted as maps $A^{\otimes k}\otimes
A \otimes A^{\otimes l}\otimes A\longrightarrow R$, which are then
denoted by the inner-product-symbol $<...>_{k,l}$ from above:
$$ <a_{1},...,a_{k+l+1},a'>_{k,l}:=(-1)^{|a'|}
   (f_{k,l}(a_{1},...,a_{k+l+1}))(a'). $$

Being an A$_\infty$-bimodule-map means that the general equation
from Proposition 4.4. is satisfied. This equation is
$$ \sum \pm f_{k,l}(...,m_{i}(...),...,a,...)+
   \sum \pm f_{k,l}(...,b_{i,j}(...,a,...),...)+ $$
$$ + \sum \pm f_{k,l}(...,a,...,m_{i}(...),...)=
   \sum \pm c_{i,j}(...,f_{k,l}(...,a,...),...). $$
Here $a\in A$ is the $(k+1)$st entry of an element in $A^{\otimes
k}\otimes A \otimes A^{\otimes l}$, which means it comes from the
\textit{A$_\infty$-bimodule} $A$, instead of the
\textit{A$_\infty$-algebra} $A$.\\
Now, by Lemma 5.1.(a), $b_{i,j}=m_{i+j+1}$ is just one of the
multiplications, and thus the left side of the equation is just
$f_{k,l}$ applied to all possible multiplications $m_{i}$. As
$f_{k,l}$ maps into $A^{*}$, one can apply the left side to an
element $a'\in A$ and therefore use the maps $<...>_{k,l}$:
\begin{equation}
  \sum \pm (f_{k,l}(...,m_{i}(...),...))(a') =
    \sum \pm <...,m_{i}(...),...,a'>_{k,l}.
\end{equation}
In order to rewrite the right side of the equation, one uses Lemma
5.1.(b) (with $f_{k,l}(...,a,...)\in A^{*}$) and the maps
$<...>_{k,l}$, when evaluating on $a'\in A$:
\begin{equation}
 \sum \pm (c_{i,j}(a_{1},...,f_{k,l}(...,a,...),...,a_{k+l+1}))(a') =
 \,\,\,\,\,\,\,\,\,\,\,\,\,\,\,\,\,\,\,\,\,\,\,\,\,\,\,\,
\end{equation}
\begin{eqnarray*}
 &=& \sum \pm (f_{k,l}(...,a,...))(m_{r}(...,a_{k+l+1},a',a_{1},...))= \\
 &=& \sum \pm <...,a,...,m_{r}(...,a_{k+l+1},a',a_{1},...)>_{k,l}.
\end{eqnarray*}

Now, with the identities (5.1) and (5.2), it is clear that the
inner-products have to satisfy equations with sums over all
possibilities of applying one multiplication to the arguments of
the inner-product subject to the conditions (i)-(iv). This is of
course just what is stated in the equation of the Proposition,
when isolating the $\partial$-terms to the left. For condition
(v), notice that the extensions of $D$ and $D^{A}$  from Lemma
2.3.(a) and Lemma 3.5.(a) record exactly the special entry $a$ in
the A$_\infty$-bimodule $A$. Thus, the A$_\infty$-bimodule element
$a$ determines the number $k$, and then $l$ is determined by the
number of arguments of $<...>_{k,l}$.

In order to see that the signs can be written as in the
proposition, one has to insert the signs for the case $m_{i}=
m_{1}=\partial$. The important terms in (5.1) are
$$ (-1)^{\varepsilon}(f_{k,l}(...,\partial a_{j},...))(a')=
   (-1)^{(\sum_{r<j} |a_{r}|)+1+k+l+1+|a'|}
   <...,\partial a_{j},...,a'>_{k,l}, $$
where $\varepsilon$ is the $\varepsilon$ from Proposition 4.4.
with $i=1$. In (5.2) one only has to look at one term, namely
$$ (-1)^{\varepsilon'} (c_{0,0}(f_{k,l}(a_{1},...,a_{k+l+1}))(a'), $$
where Proposition 4.4. implies $\varepsilon'\equiv 0 \,\,\,
(mod\,\,\, 2) $, because $i=1$. So, by Lemma 3.9., this term is
given by
\begin{eqnarray*}
 && (-1)^{|f_{k,l}(a_{1},...,a_{k+l+1})|\cdot 1}(f_{k,l}
     (a_{1},...,a_{k+l+1}))(\partial a')= \\
 && = (-1)^{k+l+(\sum_{r=1}^{k+l+1}|a_{r}|)+(|a'|-1)}
     <a_{1},...,a_{k+l+1},\partial a'>_{k,l}.
\end{eqnarray*}
Bringing this term to the left side and dividing by
$(-1)^{k+l+|a'|}$ yields the result.
\end{proof}

There is a diagrammatic way of picturing Proposition 5.4.

\begin{defn} Given an A$_\infty$-algebra $(A,D)$ with the
$\infty$-inner product $\{<.,.,...>_{k,l}:A^{\otimes k+l+2}
\longrightarrow R\}_{k\geq 0, l\geq 0}$ from Proposition 5.4. To
the inner-product $<...>_{k,l}$, one associates the symbol
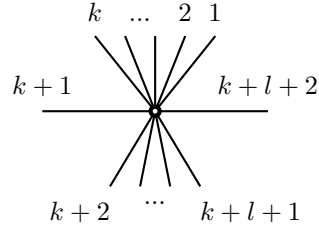
\[
\begin{pspicture}(0,0)(4,4)
 \psline(.5,2)(3.5,2)
 \psline(2,2)(1.2,3)
 \psline(2,2)(1.6,3)
 \psline(2,2)(2,3)
 \psline(2,2)(2.4,3)
 \psline(2,2)(2.8,3)
 \psline(2,2)(1.4,1)
 \psline(2,2)(1.8,1)
 \psline(2,2)(2.2,1)
 \psline(2,2)(2.6,1)
 \psdots[dotstyle=o,dotscale=2](2,2)
 \rput[b](2.8,3.2){$1$}  \rput[b](2.4,3.2){$2$}
 \rput[b](1.8,3.2){$...$} \rput[b](1.2,3.2){$k$}
 \rput[b](.5,2.2){$k+1$}
 \rput[b](3.5,2.2){$k+l+2$}
 \rput[tr](1.4,.8){$k+2$} \rput[t](2,.8){$...$} \rput[tl](2.6,.8){$k+l+1$}
\end{pspicture}
\]
More generally, to any inner-product which has (possibly iterated)
multiplications $m_{2}, m_{3}, m_{4}, ...$ (but \underline{no}
differential $\partial=m_{1}$) inside, e.g.
$$<a_{i},...,m_{j}(...),...,m_{p}(..., m_{q}(...) ,...),...>_{k,l},$$
one can associate a diagram like above, by the following rules:
\begin{itemize}
\item[i)] To every multiplication $m_{j}$, associate a tree
with $j$ inputs and one output.
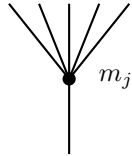
\[
\begin{pspicture}(0,0.5)(4,3.5)
 \psline(2,2)(1.2,3)
 \psline(2,2)(1.6,3)
 \psline(2,1)(2,3)
 \psline(2,2)(2.4,3)
 \psline(2,2)(2.8,3)
 \psdots[dotstyle=*,dotscale=2](2,2)
 \rput[l](2.4,2){$m_{j}$}
\end{pspicture}
\]
The symbol for the multiplication will also occur in a rotated
way. It should always be clear, where the inputs and the output
are located.
\item[ii)] To the inner product $<...>_{r,s}$, associate an "evaluation
on an open circle":
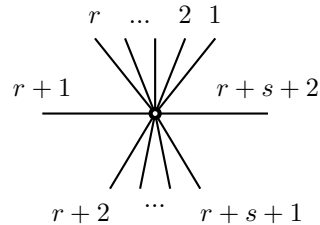
\[
\begin{pspicture}(0,0)(4,4)
 \psline(.5,2)(3.5,2)
 \psline(2,2)(1.2,3)
 \psline(2,2)(1.6,3)
 \psline(2,2)(2,3)
 \psline(2,2)(2.4,3)
 \psline(2,2)(2.8,3)
 \psline(2,2)(1.4,1)
 \psline(2,2)(1.8,1)
 \psline(2,2)(2.2,1)
 \psline(2,2)(2.6,1)
 \psdots[dotstyle=o,dotscale=2](2,2)
 \rput[b](2.8,3.2){$1$}  \rput[b](2.4,3.2){$2$}
 \rput[b](1.8,3.2){$...$} \rput[b](1.2,3.2){$r$}
 \rput[b](.5,2.2){$r+1$}
 \rput[b](3.5,2.2){$r+s+2$}
 \rput[tr](1.4,.8){$r+2$} \rput[t](2,.8){$...$} \rput[tl](2.6,.8){$r+s+1$}
\end{pspicture}
\]
Here there are $r$ elements sitting on top of the circle, $s$
elements are coming in from the bottom of the circle and the two
(special) inputs $(r+1)$ and $(r+s+2)$ on the left and right. Thus
one gets the required $r+s+2$ inputs.
\item[iii)] Around the diagram, one "sticks in" the elements
$a_{i}$ counterclockwise, (where the last element $a_{r+s+2}$ is
in the far right slot).\\
When multiplications $m_{j}$ of the graph are performed, one uses
the counterclockwise orientation of the plane to find the correct
order of the arguments $a_{i}$ in $m_{j}$ (see examples below).
\end{itemize}
Let's refer to those diagrams as
\textbf{inner-product-diagrams}.\\
Examples: Let $a$, $b$, $c$, $d$, $e$, $f$, $g$, $h$, $i$, $j$, $k$ $\in A$.\\
$<a,b,c,d>_{2,0}$, ($deg=2$):\\
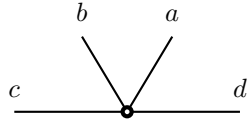
\[
\begin{pspicture}(0,0)(4,3)
 \psline(.5,1)(3.5,1)
 \psline(2,1)(1.4,2)
 \psline(2,1)(2.6,2)
 \psdots[dotstyle=o,dotscale=2](2,1)
 \rput[b](2.6,2.2){$a$}  \rput[b](1.4,2.2){$b$}
 \rput[b](0.5,1.2){$c$}  \rput[b](3.5,1.2){$d$}
\end{pspicture}
\]
$<a,b,c,d,e,f,g,h,i>_{3,4}$, ($deg=7$):\\
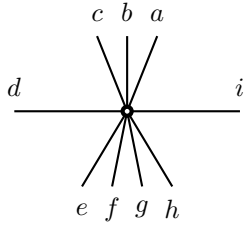
\[
\begin{pspicture}(0,0)(4,4)
 \psline(.5,2)(3.5,2)
 \psline(2,2)(1.6,3)
 \psline(2,2)(2,3)
 \psline(2,2)(2.4,3)
 \psline(2,2)(1.4,1)
 \psline(2,2)(1.8,1)
 \psline(2,2)(2.2,1)
 \psline(2,2)(2.6,1)
 \psdots[dotstyle=o,dotscale=2](2,2)
 \rput[b](2,3.2){$b$}  \rput[b](2.4,3.2){$a$}
 \rput[b](1.6,3.2){$c$}
 \rput[b](.5,2.2){$d$} \rput[b](3.5,2.2){$i$}
 \rput(1.4,.7){$e$} \rput(1.8,.7){$f$}
 \rput(2.2,.7){$g$} \rput(2.6,.7){$h$}
\end{pspicture}
\]
$<m_{2}(m_{2}(b,c),m_{2}(d,e)),m_{2}(f,a)>_{0,0}$, ($deg=0$):\\
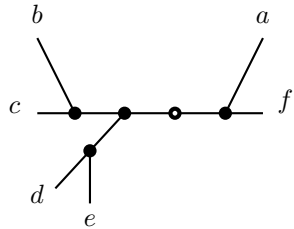
\[
\begin{pspicture}(0,0)(4,4)
 \psline(.5,2)(3.5,2)
 \psline(3,2)(3.5,3)
 \psline(1,2)(.5,3)
 \psline(1.66,2)(.74,1)
 \psline(1.2,1.5)(1.2,.8)
 \psdots[dotstyle=*,dotscale=2](1.66,2)
 \psdots[dotstyle=*,dotscale=2](1,2)
 \psdots[dotstyle=*,dotscale=2](1.2,1.5)
 \psdots[dotstyle=*,dotscale=2](3,2)
 \psdots[dotstyle=o,dotscale=2](2.33,2)
 \rput[b](.2,2){$c$}    \rput[b](3.8,2){$f$}
 \rput[b](3.5,3.2){$a$} \rput[b](0.5,3.2){$b$}
 \rput[b](.5,.8){$d$}   \rput[b](1.2,.5){$e$}
\end{pspicture}
\]
$<a,b,m_{3}(c,d,m_{2}(e,f)),g,m_{2}(h,i))>_{1,2}$, ($deg=4$):\\
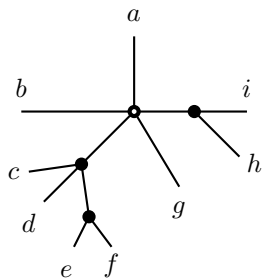
\[
\begin{pspicture}(0,-.8)(4,4)
 \psline(.5,2)(3.5,2)
 \psline(2,2)(2,3)
 \psline(2,2)(.8,.8)
 \psline(2,2)(2.6,1)
 \psline(2.8,2)(3.4,1.4)
 \psline(1.3,1.3)(.6,1.2)
 \psline(1.3,1.3)(1.4,.6)
 \psline(1.4,.6)(1.2,.2)
 \psline(1.4,.6)(1.7,.2)
 \psdots[dotstyle=o,dotscale=2](2,2)
 \psdots[dotstyle=*,dotscale=2](2.8,2)
 \psdots[dotstyle=*,dotscale=2](1.3,1.3)
 \psdots[dotstyle=*,dotscale=2](1.4,.6)
 \rput[b](2,3.2){$a$} \rput[b](.5,2.2){$b$}
 \rput(2.6,.7){$g$}   \rput[b](3.5,2.2){$i$}
 \rput[b](3.6,1.2){$h$}
 \rput[b](.4,1.1){$c$} \rput[b](.6,.4){$d$}
 \rput[b](1.1,-.2){$e$} \rput[b](1.7,-.2){$f$}
\end{pspicture}
\]
$<c,m_{2}(d,e),m_{2}(m_{2}(f,g),h),i,m_{4}(j,k,a,b)>_{2,1}$, ($deg=5$):\\
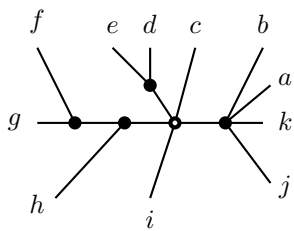
\[
\begin{pspicture}(0,0)(4,4)
 \psline(.5,2)(3.5,2)
 \psline(3,2)(3.5,3)
 \psline(1,2)(.5,3)
 \psline(1.66,2)(.74,1)
 \psline(3,2)(3.6,2.5)
 \psline(2.33,2)(2,1)
 \psline(3,2)(3.6,1.2)
 \psline(2.33,2)(2.6,3)
 \psline(2.33,2)(2,2.5)
 \psline(2,2.5)(2,3)
 \psline(2,2.5)(1.5,3)
 \psdots[dotstyle=*,dotscale=2](1.66,2)
 \psdots[dotstyle=*,dotscale=2](1,2)
 \psdots[dotstyle=*,dotscale=2](3,2)
 \psdots[dotstyle=*,dotscale=2](2,2.5)
 \psdots[dotstyle=o,dotscale=2](2.33,2)
 \rput[b](.2,1.9){$g$}  \rput[b](3.8,1.9){$k$}
 \rput[b](3.5,3.2){$b$} \rput[b](0.5,3.2){$f$}
 \rput[b](.5,.8){$h$}   \rput[b](3.8,2.5){$a$}
 \rput[b](2,.6){$i$}    \rput[b](3.8,1){$j$}
 \rput[b](2.6,3.2){$c$} \rput[b](2,3.2){$d$}
 \rput[b](1.5,3.2){$e$}
\end{pspicture}
\]
There is a chain-complex associated to the inner-product-diagrams:
\begin{itemize}
\item [i)] Degree:\\
The degree of the inner-product-diagram associated to an
inner-product $<...>_{k,l}$ with the multiplications
$m_{i_{1}},...,m_{i_{n}}$ inside is defined to be
$$ deg(Diagram):=k+l+\sum_{j=1}^{n} (i_{j}-2). $$
Examples are given above.
\item [ii)] Chain-complex:\\
For $n\geq 0$, let $C_{n}$ be the space generated by
inner-product-diagrams of degree $n$. Then let
$C:=\bigoplus_{n\geq 0} C_{n}$.
\item [iii)] Differential:\\
Let's define a differential on the inner-products to be the
composition with the operator $\tilde{\partial}:=\sum_{i}
id\otimes ...\otimes id\otimes \partial \otimes id\otimes ...
\otimes id$ (where $\partial=m_{1}$ is being in the $i$-th spot):
$$ (d(<...,m(...,m(...),...),...>))(a_{1},...,a_{s}):=
 \,\,\,\,\,\,\,\,\,\,\,\,\,\,\,\,\,\,\,\,\,\,\,\,\,\,\,
 \,\,\,\,\,\,\,\,\,\,\,\,\,\,\,\,\,\,\,\,\,\,\,\,\,\,\,$$
$$ \,\,\,\,\,\,\,\,\,\,\,\,\,\,\,\,\,\,\,\,\,\,\,\,\,\,\,
 :=(<...,m(...,m(...),...),...> ) (\sum_{i=1}^{s} (-1)^{\sum
 _{j=1}^{i-1}|a_{j}|} (a_{1},...,\partial (a_{i}),...,a_{s})). $$
Why is this well-defined and what is its diagrammatic
interpretation?\\
First let's look at the inner-product $<...>_{k,l}$ without any
multiplications inside. Then by Proposition 5.4. this means that
one puts one multiplication into the inner-product-diagram in all
possible places, such that the two lines on the far left and on
the far
right are not being multiplied (see Proposition 5.4. (iv)).\\
Now, if there are multiplications inside the inner-product, then
one can observe from Proposition 2.3., that
$$ \sum_{i} m_{n}\circ(id\otimes ... \otimes \partial\otimes ...
\otimes id) =
\,\,\,\,\,\,\,\,\,\,\,\,\,\,\,\,\,\,\,\,\,\,\,\,\,\,\,\,\,\,$$
\begin{eqnarray}
\,\,\,\,\,\,\,\,\,\,\,\,\,\,\,\,\,\,\,\,\,\,\,\,\,\,\,\,\,\,\,\,\,
\,\,\,\,\,\,\,\,\,\,\,\,\,\,\,\,\,\,\,\,\,\,\,\,\,\,\,\,\,\,\,\,\,
& = & \sum_{k=2}^{n-1} \sum_{i} m_{n+1-k} \circ(id\otimes ...
\otimes  m_{k} \otimes ... \otimes id) +\\
&& + \partial \circ m_{n}
\end{eqnarray}
(The sum over $i$ on both sides of the above equation means that
one has to put $\partial$ [or respectively $m_{k}$] in the $i$-th
spot of the tensor-product.) Now, (5.3) "brakes" the given
multiplication $m_{n}$
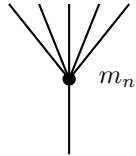
\[
\begin{pspicture}(0,0.5)(4,3.5)
 \psline(2,2)(1.2,3)
 \psline(2,2)(1.6,3)
 \psline(2,1)(2,3)
 \psline(2,2)(2.4,3)
 \psline(2,2)(2.8,3)
 \psdots[dotstyle=*,dotscale=2](2,2)
 \rput[l](2.4,2){$m_{n}$}
\end{pspicture}
\]
into all possible smaller parts $m_{n+1-k}$ and $m_{k}$
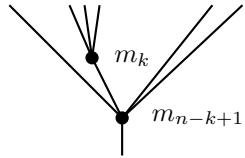
\[
\begin{pspicture}(0,0.5)(4,3.5)
 \psline(2,1)(2,1.5)
 \psline(2,1.5)(.5,3)
 \psline(2,1.5)(3.2,3)
 \psline(2,1.5)(3.6,3)
 \psline(2,1.5)(1.6,2.3)
 \psline(1.6,2.3)(1.3,3)
 \psline(1.6,2.3)(1.5,3)
 \psline(1.6,2.3)(1.7,3)
 \psdots[dotstyle=*,dotscale=2](1.6,2.3)
 \rput[l](1.9,2.3){$m_{k}$}
 \psdots[dotstyle=*,dotscale=2](2,1.5)
 \rput[l](2.4,1.5){$m_{n-k+1}$}
\end{pspicture}
\]
The last term (5.4) is also important. It makes an inductive
argument of the above possible. One gets a term $\partial
(m_{n}(...)) $ being inside the inner-product or possibly another
multiplication, which then will have arguments applied to
$\tilde{\partial}$, so that the above discussion works again.\\
So, on the level of graphs, the differential means to take just
one more multiplication in all possible places without multiplying
the given far left and far right lines (c.f. Example 5.7 below).
\end{itemize}
\end{defn}

\begin{thm} $d:C_{n}\longrightarrow C_{n-1}$, and $d^{2}=0$.
\end{thm}
\begin{proof} By the formula for the degrees in Definiton 5.5.,
a multiplication $m_{n}$ with $n$ inputs contributes by $n-2$.
Now, taking the differential means to put in one more
multiplication in all possible ways. Let's assume one wants to put
$m_{n}$ into the formula. Then this replaces $n$ arguments with
one argument in the higher level of the formula. Thus
\begin{eqnarray*}
new\,\,\,degree & = & (old\,\,\,degree)-n+1+(n-2) =\\
                & = & (old\,\,\,degree)-1.
\end{eqnarray*}

One can prove $d^{2}=0$ in two ways:
\begin{itemize}
\item[i)] Algebraically:\\
The definition of $d$ on the inner-products is just a composition
with the operator $\tilde{\partial}=\sum_{i} id\otimes ...\otimes
id\otimes \partial \otimes id\otimes ... \otimes id$ ($\partial$
being in the $i$-th spot) on $TA$. Thus $d^{2}$ is composition
with
$$ \tilde{\partial}^{2}=\sum_{i,j} \pm id\otimes ... \otimes \partial
\otimes ...\otimes \partial \otimes ... \otimes id=0. $$ This
gives zero, because $\partial$ occurring at the $i$-th and the
$j$-the spot can be obtained by first taking the one at the $i$-th
and then the one at the $j$-th, or first taking the one on the
$j$-th and then the one at the $i$-th spot. These two
possibilities cancel each other out, because $\partial$ is of
degree $-1$ and the first $\partial$ either has to "jump" over the
other $\partial$, which gives a "$-$" sign, or not.
\item[ii)] Diagrammatically (without signs):\\
If $d$ means to create one new multiplication inside the
inner-product-diagram, then $d^{2}$ obviously corresponds to
creating two new multiplications. For two given multiplications,
there are always two ways of obtaining them.
\begin{itemize}
\item[Case 1:] The multiplications are on different outputs
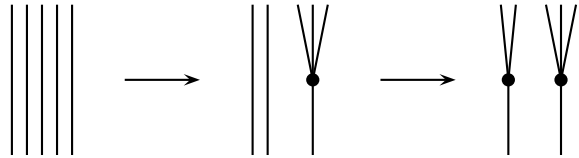
\[
\begin{pspicture}(0,0)(8.5,3)
 \psline(.5,.5)(.5,2.5)
 \psline(.7,.5)(.7,2.5)
 \psline(.9,.5)(.9,2.5)
 \psline(1.1,.5)(1.1,2.5)
 \psline(1.3,.5)(1.3,2.5)
 \psline{->}(2,1.5)(3,1.5)
 \psline(3.7,.5)(3.7,2.5)
 \psline(3.9,.5)(3.9,2.5)
 \psline(4.5,1.5)(4.3,2.5)
 \psline(4.5,.5)(4.5,2.5)
 \psline(4.5,1.5)(4.7,2.5)
 \psdots[dotstyle=*,dotscale=2](4.5,1.5)
 \psline{->}(5.4,1.5)(6.4,1.5)
 \psline(7,2.5)(7.1,1.5)
 \psline(7.2,2.5)(7.1,1.5)
 \psline(7.1,1.5)(7.1,.5)
 \psline(7.8,1.5)(7.6,2.5)
 \psline(7.8,.5)(7.8,2.5)
 \psline(7.8,1.5)(8,2.5)
 \psdots[dotstyle=*,dotscale=2](7.8,1.5)
 \psdots[dotstyle=*,dotscale=2](7.1,1.5)
\end{pspicture}
\]
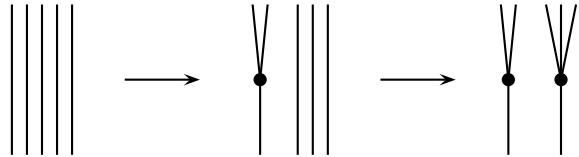
\[
\begin{pspicture}(0,0)(8.5,3)
 \psline(.5,.5)(.5,2.5)
 \psline(.7,.5)(.7,2.5)
 \psline(.9,.5)(.9,2.5)
 \psline(1.1,.5)(1.1,2.5)
 \psline(1.3,.5)(1.3,2.5)
 \psline{->}(2,1.5)(3,1.5)
 \psline(4.7,.5)(4.7,2.5)
 \psline(4.5,.5)(4.5,2.5)
 \psline(4.3,.5)(4.3,2.5)
 \psline(3.7,2.5)(3.8,1.5)
 \psline(3.9,2.5)(3.8,1.5)
 \psline(3.8,1.5)(3.8,.5)
 \psdots[dotstyle=*,dotscale=2](3.8,1.5)
 \psline{->}(5.4,1.5)(6.4,1.5)
 \psline(7,2.5)(7.1,1.5)
 \psline(7.2,2.5)(7.1,1.5)
 \psline(7.1,1.5)(7.1,.5)
 \psline(7.8,1.5)(7.6,2.5)
 \psline(7.8,.5)(7.8,2.5)
 \psline(7.8,1.5)(8,2.5)
 \psdots[dotstyle=*,dotscale=2](7.8,1.5)
 \psdots[dotstyle=*,dotscale=2](7.1,1.5)
\end{pspicture}
\]
Clearly, one can do one first and then the other, or vice versa.
In this way one always gets this term cancelling itself.
\item[Case 2:] The multiplications are on the same output\\
Here are the two ways of obtaining the same picture, and thus
cancelling out:
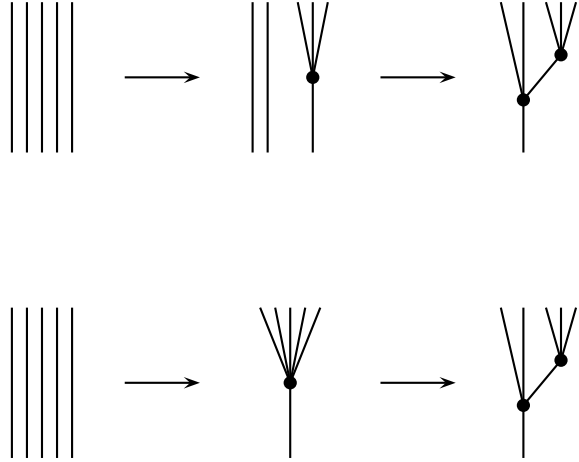
\[
\begin{pspicture}(0,0)(8.5,3)
 \psline(.5,.5)(.5,2.5)
 \psline(.7,.5)(.7,2.5)
 \psline(.9,.5)(.9,2.5)
 \psline(1.1,.5)(1.1,2.5)
 \psline(1.3,.5)(1.3,2.5)
 \psline{->}(2,1.5)(3,1.5)
 \psline(3.7,.5)(3.7,2.5)
 \psline(3.9,.5)(3.9,2.5)
 \psline(4.5,1.5)(4.3,2.5)
 \psline(4.5,.5)(4.5,2.5)
 \psline(4.5,1.5)(4.7,2.5)
 \psdots[dotstyle=*,dotscale=2](4.5,1.5)
 \psline{->}(5.4,1.5)(6.4,1.5)
 \psline(7,2.5)(7.3,1.2)
 \psline(7.3,2.5)(7.3,.5)
 \psline(7.3,1.2)(7.8,1.8)
 \psline(7.8,1.8)(7.6,2.5)
 \psline(7.8,1.8)(7.8,2.5)
 \psline(7.8,1.8)(8,2.5)
 \psdots[dotstyle=*,dotscale=2](7.8,1.8)
 \psdots[dotstyle=*,dotscale=2](7.3,1.2)
\end{pspicture}
\]
\[
\begin{pspicture}(0,0)(8.5,3)
 \psline(.5,.5)(.5,2.5)
 \psline(.7,.5)(.7,2.5)
 \psline(.9,.5)(.9,2.5)
 \psline(1.1,.5)(1.1,2.5)
 \psline(1.3,.5)(1.3,2.5)
 \psline{->}(2,1.5)(3,1.5)
 \psline(4.2,1.5)(3.8,2.5)
 \psline(4.2,1.5)(4.0,2.5)
 \psline(4.2,0.5)(4.2,2.5)
 \psline(4.2,1.5)(4.4,2.5)
 \psline(4.2,1.5)(4.6,2.5)
 \psdots[dotstyle=*,dotscale=2](4.2,1.5)
 \psline{->}(5.4,1.5)(6.4,1.5)
 \psline(7,2.5)(7.3,1.2)
 \psline(7.3,2.5)(7.3,.5)
 \psline(7.3,1.2)(7.8,1.8)
 \psline(7.8,1.8)(7.6,2.5)
 \psline(7.8,1.8)(7.8,2.5)
 \psline(7.8,1.8)(8,2.5)
 \psdots[dotstyle=*,dotscale=2](7.8,1.8)
 \psdots[dotstyle=*,dotscale=2](7.3,1.2)
\end{pspicture}
\]
\end{itemize}
\end{itemize}
\end{proof}

\begin{expl} Let $a$, $b$, $c \in A$.\\
$k=0$, $l=0$: $d(<a,b>_{0,0})=0$
\[
\begin{pspicture}(0.6,0.4)(3.6,1.8)
 \psline(1.3,1)(2.7,1)
 \psdots[dotstyle=o,dotscale=2](2,1)
 \rput[b](1.3,1.1){$a$} \rput[b](2.8,1.1){$b$}
\end{pspicture}
\]
$k=1$, $l=0$: $d(<a,b,c>_{1,0})=<a\cdot b,c>_{0,0}\pm <b,c\cdot
a>_{0,0}$
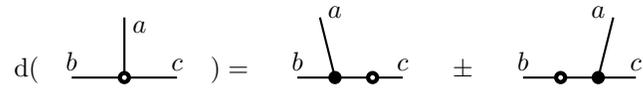
\[
\begin{pspicture}(0,0)(10,2.5)
 \psline(1.3,1)(2.7,1)
 \psline(2,1)(2,1.8)
 \psdots[dotstyle=o,dotscale=2](2,1)
 \rput[b](1.3,1.1){$b$} \rput[b](2.7,1.1){$c$} \rput[b](2.2,1.6){$a$}
 \rput(.7,1.1){$\textrm{d}($} \rput(3.4,1.1){$)=$}

 \psline(4.3,1)(5.7,1)
 \psline(4.8,1)(4.6,1.8)
 \psdots[dotstyle=*,dotscale=2](4.8,1)
 \psdots[dotstyle=o,dotscale=2](5.3,1)
 \rput[b](4.3,1.1){$b$} \rput[b](5.7,1.1){$c$} \rput[b](4.8,1.8){$a$}

 \rput(6.5,1.1){$\pm$}

 \psline(7.3,1)(8.8,1)
 \psline(8.3,1)(8.5,1.8)
 \psdots[dotstyle=*,dotscale=2](8.3,1)
 \psdots[dotstyle=o,dotscale=2](7.8,1)
 \rput[b](7.3,1.1){$b$} \rput[b](8.8,1.1){$c$} \rput[b](8.3,1.8){$a$}
\end{pspicture}
\]
$k=0$, $l=1$: $d(<a,b,c>_{0,1})=<a\cdot b,c>_{0,0}\pm <a,b\cdot
c>_{0,0}$
\[
\begin{pspicture}(0,-.2)(10,1.8)
 \psline(1.3,1)(2.7,1)
 \psline(2,1)(2,0.2)
 \psdots[dotstyle=o,dotscale=2](2,1)
 \rput[b](1.3,1.1){$a$} \rput[b](2.7,1.1){$c$} \rput[b](2.2,0.3){$b$}
 \rput(.7,1.1){$\textrm{d}($} \rput(3.4,1.1){$)=$}

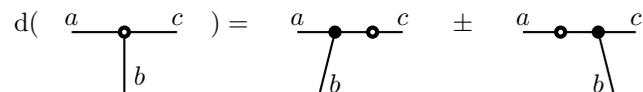
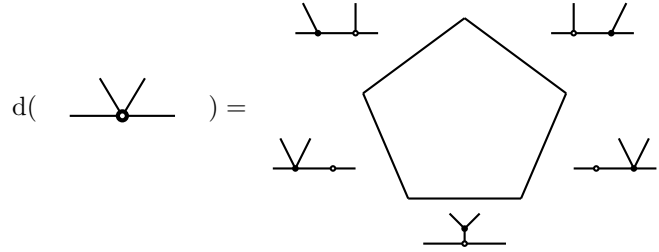
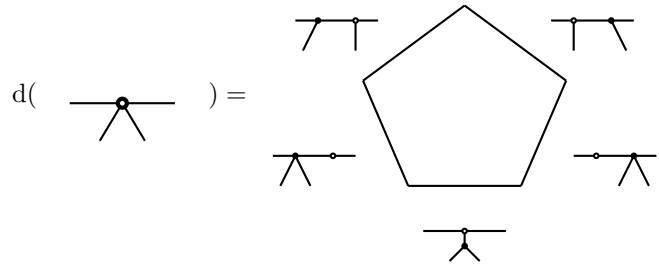
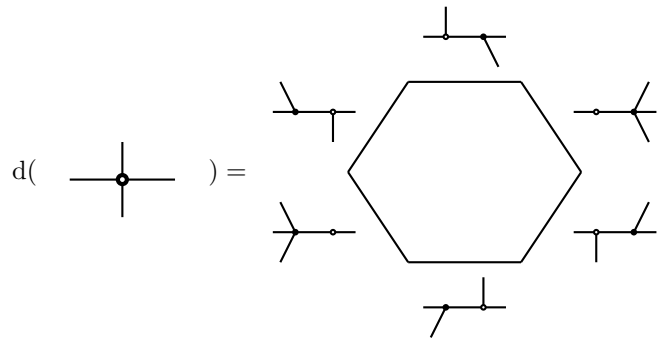
 \psline(4.3,1)(5.7,1)
 \psline(4.8,1)(4.6,0.2)
 \psdots[dotstyle=*,dotscale=2](4.8,1)
 \psdots[dotstyle=o,dotscale=2](5.3,1)
 \rput[b](4.3,1.1){$a$} \rput[b](5.7,1.1){$c$} \rput[b](4.8,0.2){$b$}

 \rput(6.5,1.1){$\pm$}

 \psline(7.3,1)(8.8,1)
 \psline(8.3,1)(8.5,0.2)
 \psdots[dotstyle=*,dotscale=2](8.3,1)
 \psdots[dotstyle=o,dotscale=2](7.8,1)
 \rput[b](7.3,1.1){$a$} \rput[b](8.8,1.1){$c$} \rput[b](8.3,0.2){$b$}
\end{pspicture}
\]
$k=2$, $l=0$:
\[
\begin{pspicture}(0,0)(10,6)
 \psline(1.3,2.9)(2.7,2.9) \psline(2,2.9)(2.3,3.4) \psline(2,2.9)(1.7,3.4)
 \psdots[dotstyle=o,dotscale=2](2,2.9)
 \rput(.7,3){$\textrm{d}($} \rput(3.4,3){$)=$}

 \psline(5.2,3.2)(6.55,4.2) \psline(6.55,4.2)(7.9,3.2)
 \psline(7.9,3.2)(7.3,1.8) \psline(7.3,1.8)(5.8,1.8)  \psline(5.8,1.8)(5.2,3.2)

 \psline(6,1.2)(7.1,1.2)
 \psline(6.55,1.4)(6.75,1.6) \psline(6.55,1.4)(6.35,1.6)
 \psdots[dotstyle=*,dotscale=1](6.55,1.4)
 \psline(6.55,1.2)(6.55,1.4) \psdots[dotstyle=o,dotscale=1](6.55,1.2)

 \psline(4,2.2)(5.1,2.2)
 \psline(4.3,2.2)(4.5,2.6) \psdots[dotstyle=*,dotscale=1](4.3,2.2)
 \psline(4.3,2.2)(4.1,2.6) \psdots[dotstyle=o,dotscale=1](4.8,2.2)

 \psline(8,2.2)(9.1,2.2)
 \psline(8.8,2.2)(8.6,2.6) \psdots[dotstyle=o,dotscale=1](8.3,2.2)
 \psline(8.8,2.2)(9,2.6) \psdots[dotstyle=*,dotscale=1](8.8,2.2)

 \psline(4.3,4)(5.4,4)
 \psline(4.6,4)(4.4,4.4) \psdots[dotstyle=*,dotscale=1](4.6,4)
 \psline(5.1,4)(5.1,4.4) \psdots[dotstyle=o,dotscale=1](5.1,4)

 \psline(7.7,4)(8.8,4)
 \psline(8,4)(8,4.4) \psdots[dotstyle=o,dotscale=1](8,4)
 \psline(8.5,4)(8.7,4.4) \psdots[dotstyle=*,dotscale=1](8.5,4)
\end{pspicture}
\]
$k=0$, $l=2$:
\[
\begin{pspicture}(0,0)(10,6)
 \psline(1.3,2.9)(2.7,2.9) \psline(2,2.9)(2.3,2.4) \psline(2,2.9)(1.7,2.4)
 \psdots[dotstyle=o,dotscale=2](2,2.9)
 \rput(.7,3){$\textrm{d}($} \rput(3.4,3){$)=$}

 \psline(5.2,3.2)(6.55,4.2) \psline(6.55,4.2)(7.9,3.2)
 \psline(7.9,3.2)(7.3,1.8) \psline(7.3,1.8)(5.8,1.8)  \psline(5.8,1.8)(5.2,3.2)

 \psline(6,1.2)(7.1,1.2)
 \psline(6.55,1)(6.75,.8) \psline(6.55,1)(6.35,.8)
 \psdots[dotstyle=*,dotscale=1](6.55,1)
 \psline(6.55,1.2)(6.55,1) \psdots[dotstyle=o,dotscale=1](6.55,1.2)

 \psline(4,2.2)(5.1,2.2)
 \psline(4.3,2.2)(4.5,1.8) \psdots[dotstyle=*,dotscale=1](4.3,2.2)
 \psline(4.3,2.2)(4.1,1.8) \psdots[dotstyle=o,dotscale=1](4.8,2.2)

 \psline(8,2.2)(9.1,2.2)
 \psline(8.8,2.2)(8.6,1.8) \psdots[dotstyle=o,dotscale=1](8.3,2.2)
 \psline(8.8,2.2)(9,1.8) \psdots[dotstyle=*,dotscale=1](8.8,2.2)

 \psline(4.3,4)(5.4,4)
 \psline(4.6,4)(4.4,3.6) \psdots[dotstyle=*,dotscale=1](4.6,4)
 \psline(5.1,4)(5.1,3.6) \psdots[dotstyle=o,dotscale=1](5.1,4)

 \psline(7.7,4)(8.8,4)
 \psline(8,4)(8,3.6) \psdots[dotstyle=o,dotscale=1](8,4)
 \psline(8.5,4)(8.7,3.6) \psdots[dotstyle=*,dotscale=1](8.5,4)
\end{pspicture}
\]
$k=1$, $l=1$:
\[
\begin{pspicture}(0,0)(10,6)
 \psline(1.3,2.9)(2.7,2.9) \psline(2,2.4)(2,3.4)
 \psdots[dotstyle=o,dotscale=2](2,2.9)
 \rput(.7,3){$\textrm{d}($} \rput(3.4,3){$)=$}

 \psline(5,3)(5.8,4.2) \psline(5.8,4.2)(7.3,4.2)  \psline(7.3,4.2)(8.1,3)
 \psline(8.1,3)(7.3,1.8) \psline(7.3,1.8)(5.8,1.8)  \psline(5.8,1.8)(5,3)

 \psline(6,1.2)(7.1,1.2)
 \psline(6.3,1.2)(6.1,.8)  \psdots[dotstyle=*,dotscale=1](6.3,1.2)
 \psline(6.8,1.2)(6.8,1.6) \psdots[dotstyle=o,dotscale=1](6.8,1.2)

 \psline(6,4.8)(7.1,4.8)
 \psline(6.3,4.8)(6.3,5.2) \psdots[dotstyle=o,dotscale=1](6.3,4.8)
 \psline(6.8,4.8)(7,4.4) \psdots[dotstyle=*,dotscale=1](6.8,4.8)

 \psline(4,2.2)(5.1,2.2)
 \psline(4.3,2.2)(4.1,2.6) \psdots[dotstyle=*,dotscale=1](4.3,2.2)
 \psline(4.3,2.2)(4.1,1.8) \psdots[dotstyle=o,dotscale=1](4.8,2.2)

 \psline(8,2.2)(9.1,2.2)
 \psline(8.3,2.2)(8.3,1.8) \psdots[dotstyle=o,dotscale=1](8.3,2.2)
 \psline(8.8,2.2)(9,2.6) \psdots[dotstyle=*,dotscale=1](8.8,2.2)

 \psline(4,3.8)(5.1,3.8)
 \psline(4.3,3.8)(4.1,4.2) \psdots[dotstyle=*,dotscale=1](4.3,3.8)
 \psline(4.8,3.8)(4.8,3.4) \psdots[dotstyle=o,dotscale=1](4.8,3.8)

 \psline(8,3.8)(9.1,3.8)
 \psline(8.8,3.8)(9,4.2) \psdots[dotstyle=o,dotscale=1](8.3,3.8)
 \psline(8.8,3.8)(9,3.4) \psdots[dotstyle=*,dotscale=1](8.8,3.8)
\end{pspicture}
\]
In the last three pictures (with $k+l=2$) the  righthand side is
understood to be a sum over the five, or respectively six, little
inner-product-diagrams. Then, as $d^{2}=0$, one can arrange these
elements according to their corresponding boundaries to give a
boundary-free object. In this way one gets certain polyhedra
associated to the $<...>_{k,l}$'s.
\end{expl}

\section{Further Remarks}

\begin{rem} The whole dicussion given here for A$_\infty$-algebras
should be transferable to different structures. It would be good
to have a notion of an $\infty$-inner-product for any algebra over
a given operad.
\end{rem}

\begin{rem}
A natural question is to ask for morphisms of A$_\infty$-algebras
preserving the $\infty$-inner-product. For this, there seem to be
several approaches. A naive way of doing this is by transferring
the usual diagrams for the condition of inner-product-preserving
maps into the "infinity-world". This seems to be too strict for
the application of compact manifolds. M. Zeinalian and the author
were able to show, that if one loosens the condition to
"inner-product-preserving maps up to homotopy", then cochains on
any compact manifold possess an $\infty$-inner-product-structure
realizing Poincar\'e-duality, which is preserved by homotopy
equivalences.\\
One way of finding a good notion of morphisms might be given by
using a general operad approach indicated in Remark 6.1.
\end{rem}

\begin{rem} There are certain well-known operations defined on the
Hochschild-cochain-complex with values in itself or its dual. Some
of these operations have a nice description for any A$_\infty
$-algebra $A$ 
(compare \cite{GJ2}).
\begin{defn} The \textbf{$\smile$-product} is a map $\smile:C^{*}(A,A)
\otimes C^{*}(A,A) \longrightarrow C^{*}(A,A)$. If $f,g \in
C^{*}(A,A)$ are given by maps $f_{j}:A^{\otimes j}\longrightarrow
A$, $g_{j}:A^{\otimes j} \longrightarrow A$, then $(f\smile g)
_{j}:A^{\otimes j}\longrightarrow A$ is given by
$$ (f\smile g)_{j}:=\sum_{k+l+m+p+q=j} \pm m_{k+m+q+2}\circ
    (id^{\otimes k} \otimes f_{l}\otimes
    id^{\otimes m} \otimes g_{p} \otimes  id^{\otimes q}). $$
\end{defn}
\begin{defn} The \textbf{Gerstenhaber-bracket} is a map $[.,.]:C^{*}(A,A)
\otimes C^{*}(A,A) \longrightarrow C^{*}(A,A)$. If $f,g \in
C^{*}(A,A)$ are given by maps $f_{j}:A^{\otimes j}\longrightarrow
A$, $g_{j}:A^{\otimes j} \longrightarrow A$, then
$$ [f,g]:=f\circ g - (-1)^{|f|\cdot|g|}g\circ f, $$
where $f\circ g$ is given by $(f\circ g) _{j}:A^{\otimes
j}\longrightarrow A$,
$$ (f\circ g)_{j}:=\sum_{k+l+m=j} \pm f_{k+1+m}\circ(id^{\otimes k}
    \otimes g_{l}\otimes id^{\otimes m}). $$
\end{defn}
Let's assume that $A$ has a unit $1\in A$.
\begin{defn} \textbf{Connes-$B$-operator} is a map $B:C^{*}(A,A^{*})
\longrightarrow C^{*}(A,A^{*})$. If $f\in C^{*}(A,A^{*})$ is given
by maps $f_{j}:A^{\otimes j}\longrightarrow A^{*}$, then $B(f)\in
C^{*}(A,A^{*})$ is given by maps $B(f)_{j-1}:A^{\otimes j-1}
\longrightarrow A^{*}$,
$$ ((B(f)_{j-1})(a_{1},...,a_{j-1}))(a_{j}):=\sum_{\sigma\in \mathbb{Z}_{j}}
    \pm (f_{j}(a_{\sigma(1)},...,a_{\sigma(j)}))(1), $$
where $\mathbb{Z}_{j}$ is understood as a subgroup of the $j$-th
symmetric group.\\
\end{defn}
The goal is now to show, that the induced $\smile$ and
$B$-operators on Hochschild-cohomology form a BV-algebra structure
whose implied Gerstenhaber structure is given by the map induced
from $[.,.]$. The author was able to proof this for (certain)
A$_\infty$-algebras with $\infty$-inner-product, by using a method
of M. Chas and D. Sullivan \cite{CS} Lemma 5.2. similar to
Gerstenhaber's proof of \cite{G} Theorem 8.5. This will be shown
in a following paper.

Ultimately, the author wants to use this to relate the
Hochschild-cohomology of a cochain-model of a compact manifold to
the homology of the free loop space of the given manifold. Here,
the BV-structure, and therefore also the Gerstenhaber-structure,
correspond to the given BV- and Gerstenhaber-structure described
in \cite{CS}. This last result was also observed in a different
way by R. Cohen and J. Jones in \cite{CJ}.
\end{rem}

Department of Mathematics, Graduate School and University Center
of the City University of New York, NY

\em E-mail address: thomastra@yahoo.com

\end{document}